\documentclass[oneside,english]{amsart}
\usepackage[T1]{fontenc}
\usepackage[utf8]{inputenc}
\usepackage{amsthm}
\usepackage{amssymb}
\usepackage{esint}
\usepackage[all]{xy}
\usepackage{lmodern}
\usepackage{color}

\makeatletter

\numberwithin{equation}{section}
\numberwithin{figure}{section}
\theoremstyle{plain}
\newtheorem{thm}{\protect\theoremname}[section]
  \theoremstyle{plain}
  \newtheorem{cor}[thm]{\protect\corollaryname}
  \newtheorem{prop}[thm]{\protect\propositionname}
  \theoremstyle{definition}
  \newtheorem{defn}[thm]{\protect\definitionname}
  \theoremstyle{remark}
  \newtheorem{rem}[thm]{\protect\remarkname}
  \theoremstyle{plain}
  
  \theoremstyle{definition}
  \newtheorem{example}[thm]{\protect\examplename}
  
  \theoremstyle{plain}

\makeatother

\usepackage{babel}
  \providecommand{\definitionname}{Definition}
  \providecommand{\examplename}{Example}
  \providecommand{\lemmaname}{Lemma}
  \providecommand{\propositionname}{Proposition}
  \providecommand{\remarkname}{Remark}
  \providecommand{\corollaryname}{Corollary}
\providecommand{\theoremname}{Theorem}

\DeclareMathOperator{\Ch}{\mathcal{C}\mathit{h}}

\begin{document}

\title{Moduli spaces of (bi)algebra structures in topology and geometry}

\author{Sinan Yalin}

\email{yalinprop@gmail.com}

\begin{abstract}
After introducing some motivations for this survey, we describe a formalism to parametrize a wide class of algebraic structures occurring naturally in various problems of topology, geometry and mathematical physics. This allows us to define an ``up to homotopy version'' of algebraic structures which is coherent (in the sense of $\infty$-category theory) at a high level of generality. To understand the classification and deformation theory of these structures on a given object, a relevant idea inspired by geometry is to gather them in a moduli space with nice homotopical and geometric properties. Derived geometry provides the appropriate framework to describe moduli spaces classifying objects up to weak equivalences and encoding in a geometrically meaningful way their deformation and obstruction theory. As an instance of the power of such methods, I will describe several results of a joint work with Gregory Ginot related to longstanding conjectures in deformation theory of bialgebras, $E_n$-algebras and quantum group theory.
\end{abstract}

\maketitle

\tableofcontents{}

\section*{Introduction}

To motivate a bit the study of algebraic structures and their moduli spaces in topology, we will simply start from singular cohomology. Singular cohomology provides a first approximation of the topology of a given space by its singular simplices, nicely packed in a cochain complex. Computing the cohomology of spaces already gives us a way to distinguish them and extract some further information like characteristic classes for instance. Singular cohomology has the nice property to be equipped with an explicit commutative ring structure given by the cup product. This additional structure can distinguish spaces which have the same cohomology groups, illustrating of the following idea: adding finer algebraic structures is a way to parametrize finer invariants of our spaces. In the case of manifolds, it can also be used to get more geometric data (from characteristic classes and Poincaré duality for instance).
Such an algebraic struture determined by operations with several inputs and one single output (the cup product in our example) satisfying relations (associativity, commutativity) is parametrized by an operad (here the operad $Com$ of commutative associative algebras). More generally, the notion of operad has proven to be a fundamental tool to study algebras playing a key role in algebra, topology, category theory, differential and algebraic geometry, mathematical physics (like Lie algebras, Poisson algebras and their variants).

We can go one step further and relax such structures \emph{up to homotopy} in an appropriate sense. Historical examples for this include higher Massey products, Steenrod squares and (iterated) loop spaces.

Higher Massey products organize into an $A_{\infty}$-algebra structure on the cohomology of a space and give finer invariants than the cup product. For instance, the trivial link with three components has the same cohomology ring as the borromean link (in both cases, the cup product is zero), but the triple Massey product vanishes in the second case and not in the first one, implying these links are not equivalent.

Loop spaces are another fundamental example of $A_{\infty}$-algebras (in topological spaces this time). When one iterates this construction by taking the loop space of the loop space and so on, one gets an $E_n$-algebra (more precisely an algebra over the little $n$-disks operad). These algebras form a hierarchy of ``more and more'' commutative and homotopy associative structures, interpolating between $A_\infty$-algebras (the $E_1$ case, encoding homotopy associative structures) and $E_{\infty}$-algebras (the colimit of the $E_n$'s, encoding homotopy commutative structures). Algebras governed by $E_n$-operads and their deformation theory play a prominent role in a variety of topics, not only the study of iterated loop spaces but also Goodwillie-Weiss calculus for embedding spaces, deformation quantization of Poisson manifolds, Lie bialgebras and shifted Poisson structures in derived geometry,and factorization homology of manifolds~\cite{Ko1, Ko2, Lur0, Lur2, CPTVV, FG, Fra, Fre5, GJ, GTZ, Hin, Kap-TFT, May, Preygel, Shoikhet, Tam1, Toen-ICM}.

The cup product is already defined at a chain level but commutative only up to homotopy, meaning that there is an infinite sequence of obstructions to commutativity given by the so called higher cup products. That is, these higher cup products form an $E_{\infty}$-algebra structure on the singular cochain complex. This $E_{\infty}$-structure classifies the rational homotopy type of spaces (this comes from Sullivan's approach to rational homotopy theory \cite{Sul}) and the integral homotopy type of finite type nilpotent spaces (as proved by Mandell \cite{Man}). Moreover, such a structure induces the Steenrod squares acting on cohomology and, for Poincaré duality spaces like compact oriented manifolds for example, the characteristic classes that represent these squares (the Wu classes). This is a first instance of how a homotopy algebraic structure can be used to build characteristic classes. Then, to study operations on generalized cohomology theories, one moves from spaces to the stable homotopy theory of spectra, the natural recipient for (generalized) cohomology theories, and focuses on theories represented by $E_{\infty}$-ring spectra (more generally, highly structured ring spectra). In this setting, the moduli space approach (in a homotopy theoretic way) already proved to be useful \cite{GoH} (leading to a non trivial improvement of the Hopkins-Miller theorem in the study of highly structured ring spectra).

However, algebraic structures not only with products but also with coproducts, play a crucial role in various places in topology, geometry and mathematical physics. One could mention for instance the following important examples: Hopf algebras in representation theory and mathematical physics, Frobenius algebras encompassing the Poincaré duality phenomenon in algebraic topology and deeply related to field theories, Lie bialgebras introduced by Drinfeld in quantum group theory, involutive Lie bialgebras as geometric operations on the equivariant homology of free loop spaces in string topology. A convenient way to handle such kind of structures is to use the formalism of props, a generalization of operads encoding algebraic structures based on operations with several inputs and several outputs.

A natural question is then to classify such structures (do they exist, how many equivalence classes) and to understand their deformation theory (existence of infinitesimal perturbations, formal perturbations, how to classify the possible deformations). Understanding how they are rigid or how they can be deformed provide informations about the objects on which they act and new invariants for these objects. For this, a relevant idea inspired from geometry come to the mind, the notion of moduli space, a particularly famous example being the moduli spaces of algebraic curves (or Riemann surfaces). The idea is to associate, to a collection of objects we want to parametrize equipped with an equivalence relation (surfaces up to diffeomorphism, vector bundles up to isomorphism...), a space $\mathcal{M}$ whose points are these objects and whose connected components are the equivalence classes of such objects. This construction is also called a classifying space in topology, a classical example being the classifying space $BG$ of a group $G$, which parametrizes isomorphism classes of principal $G$-bundles. If we are interested also in the deformation theory of our collection of objects (how do we allow our objects to modulate), we need an additional geometric structure which tells us how we can move infinitesimally our points (tangent spaces). To sum up, the guiding lines of the moduli space approach are the following:
\begin{itemize}
\item[$\bullet$] To determine the non-emptiness of $\mathcal{M}$ and to compute $\pi_*\mathcal{M}$ solve existence and unicity problems;

\item[$\bullet$] The geometric structure of $\mathcal{M}$ imply the existence of tangent spaces. The tangent space over a given point $x$ of $\mathcal{M}$ is a dg Lie algebra controlling the (derived) deformation theory of $x$ (deformations of $x$ form a derived moduli problem) in a sense we will precise in Section 3;

\item[$\bullet$] One can ``integrate'' over $\mathcal{M}$ to produce invariants of the objects parametrized by $\mathcal{M}$. Here, the word ``integrate'' has to be understood in the appropriate sense depending on the context: integrating a differential form, pairing a certain class along the (virtual) fundamental class of $\mathcal{M}$, etc.
\end{itemize}
We already mentioned \cite{GH} (inspired by the method of \cite{BDG}) as an application of the first item in the list above. In the second one, we mention \emph{derived} deformation theory and \emph{derived} moduli problems, which implicitely assume that, in some sense, our moduli space $\mathcal{M}$ lives in a ($\infty$-)category of derived objects (e.g. derived schemes, derived stacks...) where tangent spaces are actually complexes. This is due to the fact that we want to encompass the whole deformation theory of points, and this cannot be done in the classical setting: for varieties or schemes, the tangent space is a vector space which consists just of the equivalence classes of infinitesimal deformations of the point. For stacks, the tangent space is a two-term complex whose $H^0$ is the set of equivalences classes of infinitesimal deformations, and $H^{-1}$ is the Lie algebra of automorphisms of the point (infinitesimal automorphisms). But obstruction theory does not appear on the tangent structure here, because it has to live in positive degrees (we will go back to this remark in Section 3).

As a last remark about the third item in the list above, let us say that the idea of using moduli spaces to produce topological invariants got also a lot of inspiration from quantum field theory and string theory in the 80's. By the Feynman path integral approach, the equations describing the evolution of a quantum system are determined by the minimas of a functional integral over all the possible paths of this system, that is, by integrating a certain functional over the space of fields. This independence from a choice of path led to the idea that one could build a topological invariant of a geometric object by computing an integral over the moduli space of all possible geometric structures of this kind, ensuring automatically the desired invariance property. This is the principle underlying two important sorts of invariants of manifolds. First, Witten's quantization of the classical Chern-Simmons invariant in the late 80s \cite{Wit}, which provided topological invariants for $3$-dimensional manifolds (including known invariants such as the Jones polynomial) by integrating a geometric invariant over a moduli space of connections. Second, Kontsevich's formalization of Gromov– Witten invariants in symplectic topology (counting pseudo-holomorphic curves) and algebraic geometry (counting algebraic curves), defined by a pairing along the virtual fundamental class of the moduli space (stack) of stable maps (an analogue of the fundamental class for singular objects suitably embedded in a derived setting).

\textbf{Organization of the paper.} The first section is devoted to the formalism of props and algebras over props, accompanied by relevant examples of topological or geometric origin in Section 2. The third section focuses on algebraic structures up to homotopy, defined as algebras over a cofibrant resolution of the prop, and the fundamental theorem asserting that this notion does not rely on the choice of such a resolution (up to an equivalence of $\infty$-categories). This lays down the coherent foundations to study homotopy bialgebras. We then provide a little introduction to derived algebraic geometry and formal moduli problems in Section 4, without going too far in the details (we refer the reader to \cite{Toe} for a more thorough survey on this topic), before formalizing the idea of moduli spaces of algebraic structures in Section 5 as well as their most important properties. The way to recover geometrically deformation theory and obstruction theory for such structures is explained in Section 5. Section 6 describes a joint work with Gregory Ginot, merging the homotopical and geometric theory of such moduli spaces with several features of factorization homology (higher Hochschild (co)homology) to solve several open conjectures in deformation theory of $E_n$-algebras and bialgebras related to quantum group theory.

\textbf{Acknowledgements}. The idea of writing such a survey originates in the inaugural two-weeks program at the mathematical research institute MATRIX in Australia called Higher Structures in Geometry and Physics, which took place in June 2016. The author gave a talk at this program about moduli spaces of algebraic structures and their application to the recent paper \cite{GY}. The present article is somehow a (largely) extended version of his talk, which will be eventually part of a Proceedings Volume devoted to this workshop.  The author would like to thank the MATRIX institute for supporting this program, the organizers of this programme for inviting him, and all the participants for their interest and for the very enjoyable atmosphere during the two weeks spent there. Last but not least, kangaroos are very much thanked for their natural awesomeness.

\section{Parametrizing algebraic structures}

To simplify our exposition, we will work in the base category $\Ch$ of $\mathbb{Z}$-graded cochain complexes over a field $\mathbb{K}$ of characteristic zero. Before stating the general definition of a prop, let us give a few examples of algebraic structures the reader may have encountered already.
\begin{example}
Differential graded (dg for short) associative algebras are complexes $A$ equipped with an associative product $A\otimes A\rightarrow A$. We can represent such an operation by an oriented graph with two inputs and one output $\xymatrix @R=0.5em@C=0.75em@M=0em{
\ar@{-}[dr] & & \ar@{-}[dl] \\
 & \ar@{-}[d] &  \\
 & & }$
satisfying the associativity relation
\[
\xymatrix @R=0.5em@C=0.75em@M=0em{
 \ar@{-}[dr] & & \ar@{-}[dr] & & \ar@{-}[dl] & &  \ar@{-}[dr] & & \ar@{-}[dl] & & \ar@{-}[dl]\\
 & \ar@{-}[dr] & &\ar@{-}[dl] & & = &  & \ar@{-}[dr] & &\ar@{-}[dl] & \\
 & & \ar@{-}[d] & & & & & & \ar@{-}[d] & & & \\
 & & & & & & & & & & \\
}
\]
Common examples of such structures include algebras $\mathbb{K}[G]$ of finite groups $G$ in representation theory, or the singular cochains $C^*(X;\mathbb{Z})$ of a topological space equipped with the cup product $\cup$ of singular simplices. In the first case we have an associative algebra in $\mathbb{K}$-modules, in the second case this is a dg associative algebra, so the cup product is a cochain morphism determined by linear maps
\[
\cup:C^m(X;\mathbb{Z})\otimes C^n(X;\mathbb{Z})\rightarrow C^{m+n}(X;\mathbb{Z}).
\]
\end{example}
\begin{example}
In certain cases, the product is not only associative but also commutative, and one call such algebras commutative dg algebras or cdgas. To represent graphically this symmetry condition, we index the inputs of the product
$\xymatrix @R=0.5em@C=0.75em@M=0em{
1\ar@{-}[dr] & & 2\ar@{-}[dl]  \\
 & \ar@{-}[d] &  \\
  & 1 & }$
and add the symmetry condition
\[
\xymatrix @R=0.5em@C=0.75em@M=0em{
1\ar@{-}[dr] & & 2\ar@{-}[dl]  \\
 & \ar@{-}[d] &  \\
  & 1 & }=
\xymatrix @R=0.5em@C=0.75em@M=0em{
2\ar@{-}[dr] & & 1\ar@{-}[dl]  \\
 & \ar@{-}[d] &  \\
  & 1 & }.
  \]
A way to rephrase this symmetry is to say that $\Sigma_2$ acts trivially on $\xymatrix @R=0.5em@C=0.75em@M=0em{
1\ar@{-}[dr] & & 2\ar@{-}[dl]  \\
 & \ar@{-}[d] &  \\
  & 1 & }$.
In the dg setting, this symmetry has to be understood in the graded sense, that is $ab=(-1)^{deg(a)deg(b)}ba$.
Commutative algebras are very common objects, for instance the singular cohomology of spaces equipped with the cup product defined previously at the chain level, or the de Rham cohomology for manifolds. Commutative rings also represent affine schemes in algebraic geometry or rings of functions on differentiable manifolds. Cdgas over $\mathbb{Q}$ also model the rational homotopy type of simply connected spaces.
\end{example}
\begin{example}
Another example of product is the bracket defining Lie algebras
$\xymatrix @R=0.5em@C=0.75em@M=0em{
1\ar@{-}[dr] & & 2\ar@{-}[dl]  \\
 & \ar@{-}[d] &  \\
  & 1 & }$
satisfying an antisymmetry condition
\[
\xymatrix @R=0.5em@C=0.75em@M=0em{
1\ar@{-}[dr] & & 2\ar@{-}[dl]  \\
 & \ar@{-}[d] &  \\
  & 1 & }= -
\xymatrix @R=0.5em@C=0.75em@M=0em{
2\ar@{-}[dr] & & 1\ar@{-}[dl]  \\
 & \ar@{-}[d] &  \\
  & 1 & }
  \]
and the Jacobi identity
\[
\xymatrix @R=0.5em@C=0.75em@M=0em{
 1\ar@{-}[dr] & & 2\ar@{-}[dr] & & 3\ar@{-}[dl]\\
 & \ar@{-}[dr] & &\ar@{-}[dl] &  \\
 & & \ar@{-}[d] & & \\
 & & & & \\
}
+
\xymatrix @R=0.5em@C=0.75em@M=0em{
 3\ar@{-}[dr] & & 1\ar@{-}[dr] & & 2\ar@{-}[dl]\\
 & \ar@{-}[dr] & &\ar@{-}[dl] &  \\
 & & \ar@{-}[d] & & \\
 & & & & \\
}
+
\xymatrix @R=0.5em@C=0.75em@M=0em{
 2\ar@{-}[dr] & & 3\ar@{-}[dr] & & 1\ar@{-}[dl]\\
 & \ar@{-}[dr] & &\ar@{-}[dl] &  \\
 & & \ar@{-}[d] & & \\
 & & & & \\
}=0.
\]
A way to rephrase the antisymmetry is to say that the action of $\Sigma_2$ on $\xymatrix @R=0.5em@C=0.75em@M=0em{
1\ar@{-}[dr] & & 2\ar@{-}[dl]  \\
 & \ar@{-}[d] &  \\
  & 1 & }$ is given by the signature representation $sgn_2$. Lie algebras appear for instance as tangent spaces of Lie groups in differential geometry (Lie's third theorem gives an equivalence between the category of finite dimensional Lie algebras in vector spaces and the category of simply connected Lie groups), in Quillen's approach to rational homotopy theory and in deformation theory (``Deligne principle'' relating formal moduli problems to dg Lie algebras).
\end{example}
In these three first examples, we see that the algebraic structure is defined only by operations with several inputs and one single output. Such structures can be encoded by a combinatorial object called an operad, and a given kind of algebra is an algebra over the associated operad. We refer the reader to \cite{LV} for more details about this formalism. However, there are more general algebraic structures involving operations with several inputs and several outputs. We give below two fundamental examples of these, before unwrapping the general definition of the combinatorial structure underlying them (props).
\begin{example}
Poisson-Lie groups are Lie groups with a compatible Poisson structure, which occur in mathematical physics as gauge groups of certain classical mechanical systems such as integrable systems. Because of the Poisson bracket, the tangent space $T_eG$ of a Poisson-Lie group $G$ at the neutral element $e$ is equipped with  a ``Lie cobracket'' compatible with its Lie algebra structure, so that $T_eG$ forms something called a Lie bialgebra.
The compatibility relation between the bracket and the cobracket is called the Drinfeld's compatibility relation or the cocycle relation. In terms of graphical presentation, we have a bracket and a cobracket
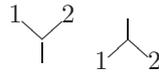
\[
\xymatrix @R=0.5em@C=0.75em@M=0em{
1\ar@{-}[dr] & & 2\ar@{-}[dl] & & \ar@{-}[d] & \\
 & \ar@{-}[d] & & & \ar@{-}[dl] \ar@{-}[dr] & \\
 & & & 1 & & 2 }
\]
which are antisymmetric, that is, with the signature action of $\Sigma_2$.
These two operations satisfy the following relations:
\begin{itemize}
\item[]Jacobi
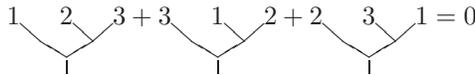
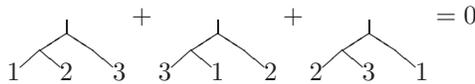
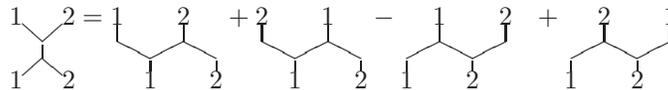
\[
\xymatrix @R=0.5em@C=0.75em@M=0em{
 1\ar@{-}[dr] & & 2\ar@{-}[dr] & & 3\ar@{-}[dl]\\
 & \ar@{-}[dr] & &\ar@{-}[dl] &  \\
 & & \ar@{-}[d] & & \\
 & & & & \\
}
+
\xymatrix @R=0.5em@C=0.75em@M=0em{
 3\ar@{-}[dr] & & 1\ar@{-}[dr] & & 2\ar@{-}[dl]\\
 & \ar@{-}[dr] & &\ar@{-}[dl] &  \\
 & & \ar@{-}[d] & & \\
 & & & & \\
}
+
\xymatrix @R=0.5em@C=0.75em@M=0em{
 2\ar@{-}[dr] & & 3\ar@{-}[dr] & & 1\ar@{-}[dl]\\
 & \ar@{-}[dr] & &\ar@{-}[dl] &  \\
 & & \ar@{-}[d] & & \\
 & & & & \\
}=0
\]
\item[]co-Jacobi
\[
\xymatrix @R=0.5em@C=0.75em@M=0em{
 & & \ar@{-}[d] & & \\
 & &\ar@{-}[dl] \ar@{-}[dr]& & \\
 &\ar@{-}[dl] \ar@{-}[dr] & &\ar@{-}[dr] & \\
 1 & & 2 & & 3  \\
}
+
\xymatrix @R=0.5em@C=0.75em@M=0em{
 & & \ar@{-}[d] & & \\
 & &\ar@{-}[dl] \ar@{-}[dr]& & \\
 &\ar@{-}[dl] \ar@{-}[dr] & &\ar@{-}[dr] & \\
 3 & & 1 & & 2 \\
}
+
\xymatrix @R=0.5em@C=0.75em@M=0em{
 & & \ar@{-}[d] & & \\
 & &\ar@{-}[dl] \ar@{-}[dr]& & \\
 &\ar@{-}[dl] \ar@{-}[dr] & &\ar@{-}[dr] & \\
 2 & & 3 & & 1 \\
}=0
\]
\item[]The cocycle relation
\[
\xymatrix @R=0.5em@C=0.75em@M=0em{
1\ar@{-}[dr] & & 2\ar@{-}[dl]  \\
 & \ar@{-}[d] & \\
 & \ar@{-}[dl] \ar@{-}[dr] &  \\
 1 & & 2 \\
}
=
\xymatrix @R=0.5em@C=0.75em@M=0em{
 1\ar@{-}[d] & & 2\ar@{-}[d] & \\
 \ar@{-}[dr] & & \ar@{-}[dl] \ar@{-}[dr] & \\
 & \ar@{-}[d] & & \ar@{-}[d] \\
 & 1 & & 2 \\
}
+
\xymatrix @R=0.5em@C=0.75em@M=0em{
 2\ar@{-}[d] & & 1\ar@{-}[d] & \\
 \ar@{-}[dr] & & \ar@{-}[dl] \ar@{-}[dr] & \\
 & \ar@{-}[d] & & \ar@{-}[d] \\
 & 1 & & 2 \\
}
-
\xymatrix @R=0.5em@C=0.75em@M=0em{
 &1\ar@{-}[d] & & 2\ar@{-}[d] \\
 & \ar@{-}[dl] \ar@{-}[dr]& & \ar@{-}[dl] \\
 \ar@{-}[d] & & \ar@{-}[d] & & \\
 1& &2 & \\
}
+
\xymatrix @R=0.5em@C=0.75em@M=0em{
 &2\ar@{-}[d] & & 1\ar@{-}[d] \\
 & \ar@{-}[dl] \ar@{-}[dr]& & \ar@{-}[dl] \\
 \ar@{-}[d] & & \ar@{-}[d] & & \\
1 & &2 & \\
}
\]
\end{itemize}
The cocycle relation means that the Lie cobracket of a Lie bialgebra $g$ is a cocycle in the Chevalley-Eilenberg
complex $C^*_{CE}(g,\Lambda^2g)$, where $\Lambda^2g$ is equipped with the structure of
$g$-module induced by the adjoint action. Let us note that there is an analogue of Lie's third theorem in this context, namely the category of finite dimensional Lie bialgebras in vector spaces is equivalent to the category of simply connected Poisson-Lie groups \cite{Dri2}. Deformation quantization of Lie bialgebras produces quantum groups, which turned out to be relevant for mathematical physics and for low-dimensional topology (quantum invariants of knots and $3$-manifolds). This process also deeply involves other kind of objects such as Grothendieck-Teichmüller groups, multizeta values via the Drinfeld associators \cite{Dri2}, or graph complexes. The problem of a universal quantization of Lie bialgebras raised by Drinfeld was solved by Etingof and Kazhdan \cite{EK1}, \cite{EK2}. A deformation quantization of a Lie bialgebra $g$ is a topologically free Hopf algebra $H$ over the ring of formal power series $\mathbb{K}[[\hbar]]$ such that $H/\hbar H$ is isomorphic to $U(g)$ (the enveloping algebra of $g$)as a co-Poisson bialgebra. Such a Hopf algebra is called a quantum universal enveloping algebra (QUE for short). The general idea underlying this process is to tensor the $\mathbb{K}$-linear category of $g$-modules by formal power series, equip it with a braided monoidal structure induced by the choice of a Drinfeld associator and an $r$-matrix, and make the forgetful functor from $g[[\hbar]]$-modules to $\mathbb{K}[[\hbar]]$-modules braided monoidal. Applying the tannakian formalism to this functor, the category of $g[[\hbar]]$-modules is equivalent to the category of modules over the QUE algebra of $g$. Deformation quantization of Lie bialgebras can be formulated in the formalism of props and their algebras, see for instance the introduction of \cite{EE} explaining quantization/de-quantization problems in terms of prop morphisms. Another point of view is the prop profile approach of \cite{Mer2}, particularly useful to relate the results of \cite{GY} to deformation quantization of Lie bialgebras.

A variant of Lie bialgebras called involutive Lie bialgebras arose in low dimensional topology, in the work of Goldman \cite{Gol2} and Turaev \cite{Tur}. Given a surface $S$, one considers the $\mathbb{K}$-module generated by the free homotopy classes of loops on $S$. Let us note $L:S^1\rightarrow S$ a free loop on the surface $S$ (that is, a continuous map which is not pointed, contrary to based loops) and $[L]$ its free homotopy class. Up to homotopy, we can make two loops intersect transversely, so we suppose that two given loops $L$ and $K$ intersect only at a finite number of points, and we note $L\cap K$ this finite set. The Lie bracket of $[L]$ and $[K]$ is then defined by
\[
\{[L],[K]\} =\sum_{p\in L\cap K}\epsilon_p [L\cup_p K]
\]
where $L\cup_p K$ is the loop parametrized by going from $p$ to $p$ along $L$, then going again from $p$ to $p$ along $K$. The symbol $\epsilon_p$ denote a number which is $-1$ or $1$, depending on the way $L$ and $K$ intersect at $p$ with respect to a choice of orientation. The cobracket is then defined similarly, by considering this time the self-intersections of $L$ (that we can take transverse, up to homotopy):
\[
\delta ([L])=\sum_{p\in L\cap L}\epsilon_p ([L_{1,p}]\otimes [L_{2,p}] - [L_{2,p}]\otimes [L_{1,p}])
\]
where $L_1$ and $L_2$ are the two loops obtained by separating $L$ in two parts at the self-intersection point $p$.
These two operations define a Lie bialgebra structure, satisfying moreover $(\{,\}\circ\delta)([L])=0$. From the graphical presentation viewpoint, this means that we add the involutivity relation
\[
\xymatrix @R=0.5em@C=0.75em@M=0em{
 & \ar@{-}[d] & \\
 & \ar@{-}[dl] \ar@{-}[dr] & \\
 \ar@{-}[dr] & & \ar@{-}[dl] \\
 & \ar@{-}[d] & \\
 & & \\
}=0.
\]
Links defined in the cylinder $S\times [0;1]$ over $S$ can be presented by diagrams of loops on $S$ via the canonical projection $S\times [0;1]\rightarrow S$. Explicit quantizations of the Lie bialgebra of loops on $S$ have been used to produce (quantum) invariants of those links \cite{Tur} and the corresponding $3$-dimensional TQFTs \cite{RT}.

Ten years after, algebraic structures on free loop spaces for more general manifolds were introduced by Chas and Sullivan, giving birth to string topology \cite{CS1}, a very active field of research nowadays. In the equivariant setting, the Lie bialgebra of Goldman and Turaev has been generalized to loop spaces of smooth manifolds \cite{CS2}.
The string homology of a smooth manifold  $M$ is defined as the reduced equivariant homology (i.e. relative to constant loops) of the free loop space $LM$ of $M$. The word equivariant refers here to the action of $S^1$ on loops
by rotation. According to \cite{CS2}, the string homology of a smooth manifold forms an involutive Lie
bialgebra. Let us note that for an $n$-dimensional manifold, the bracket and the cobracket of this
structure are of degree $2-n$. In particular, the string homology of a surface is isomorphic to Goldman-Turaev Lie bialgebra as a graded Lie bialgebra. Let us note that such a structure is also related to very active research topics in symplectic topology. Precisely, the string homology of $M$ is isomorphic as a graded Lie bialgebra to the contact homology of its cotangent bundle (equipped with the standard symplectic form) \cite{CL}. This result is part of a larger program aimed at relating string topology and symplectic field theory.
\end{example}
\begin{example}
A dg Frobenius algebra is a unitary dg commutative associative algebra of finite dimension $A$ endowed with a symmetric non-degenerated bilinear form $<.,.>:A\otimes A\rightarrow \mathbb{K}$
which is invariant with respect to the product, i.e $<xy,z>=<x,yz>$.

A dg Frobenius bialgebra of degree $m$ is a triple $(B,\mu,\Delta)$ such that:

(i) $(B,\mu)$ is a dg commutative associative algebra;

(ii) $(B,\Delta)$ is a dg cocommutative coassociative coalgebra with $deg(\Delta)=m$;

(iii) the map $\Delta:B\rightarrow B\otimes B$ is a morphism of left $B$-module
and right $B$-module, i.e in Sweedler's notations we have the Frobenius relations
\begin{eqnarray*}
\sum_{(x.y)}(x.y)_{(1)}\otimes (x.y)_{(2)} & = & \sum_{(y)}x.y_{(1)}\otimes y_{(2)}\\
 & = & \sum_{(x)}(_1)^{m|x|}x_{(1)}\otimes x_{(2)}.y
\end{eqnarray*}

The two definitions are strongly related. Indeed, if $A$ is a Frobenius algebra, then the pairing
$<.,.>$ induces an isomorphism of $A$-modules $A\cong A^*$, hence a map
\[
\Delta:A\stackrel{\cong}{\rightarrow} A^* \stackrel{\mu^*}{\rightarrow} (A\otimes A)^*\cong A^*\otimes
A^*\cong A\otimes A
\]
which equips $A$ with a structure of Frobenius bialgebra.
Conversely, one can prove that every unitary counitary Frobenius bialgebra gives rise to a Frobenius algebra,
so the two notions are equivalent. In terms of graphical presentation, we have a product of degree $0$ and a coproduct of degree $m$ presented by
\[
\xymatrix @R=0.5em@C=0.75em@M=0em{
\ar@{-}[dr] & & \ar@{-}[dl] & & \ar@{-}[d] & \\
 & \ar@{-}[d] & & & \ar@{-}[dl] \ar@{-}[dr] & \\
 & & & & & }
\]
and satisfying the following relations:
\begin{itemize}
\item[]Associativity and coassociativity
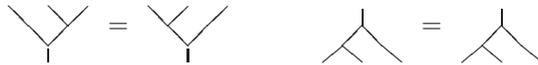
\[
\xymatrix @R=0.5em@C=0.75em@M=0em{
 \ar@{-}[dr] & & \ar@{-}[dr] & & \ar@{-}[dl] & &  \ar@{-}[dr] & & \ar@{-}[dl] & & \ar@{-}[dl]\\
 & \ar@{-}[dr] & &\ar@{-}[dl] & & = &  & \ar@{-}[dr] & &\ar@{-}[dl] & \\
 & & \ar@{-}[d] & & & & & & \ar@{-}[d] & & & \\
 & & & & & & & & & & \\
}
\hspace*{1cm}
\xymatrix @R=0.5em@C=0.75em@M=0em{
 & & \ar@{-}[d] & & & & & & \ar@{-}[d] & & & \\
 & &\ar@{-}[dl] \ar@{-}[dr]& & & = &  & &\ar@{-}[dl] \ar@{-}[dr] & & \\
 &\ar@{-}[dl] \ar@{-}[dr] & &\ar@{-}[dr] & & & & \ar@{-}[dr] \ar@{-}[dl]& & \ar@{-}[dr] & \\
 & & & & & & & & & & \\
}
\]

\item[]Frobenius relations
\[
\xymatrix @R=0.5em@C=0.75em@M=0em{
\ar@{-}[dr] & & \ar@{-}[dl] & & \ar@{-}[d] & & \ar@{-}[d] & \\
 & \ar@{-}[d] & & = & \ar@{-}[dr] & & \ar@{-}[dl] \ar@{-}[dr] & \\
 & \ar@{-}[dl] \ar@{-}[dr] & & & & \ar@{-}[d] & & \ar@{-}[d] \\
 & & & & & & & \\
}
\hspace*{1cm}
\xymatrix @R=0.5em@C=0.75em@M=0em{
\ar@{-}[dr] & & \ar@{-}[dl] & & &\ar@{-}[d] & & \ar@{-}[d] \\
 & \ar@{-}[d] & & = & & \ar@{-}[dl] \ar@{-}[dr]& & \ar@{-}[dl] \\
 & \ar@{-}[dl] \ar@{-}[dr] & & & \ar@{-}[d] & & \ar@{-}[d] & \\
 & & & & & & & \\
}
\]
\end{itemize}
In the unitary and counitary case, one adds a generator for the unit, a generator for the counit and the necessary
compatibility relations with the product and the coproduct. We refer the reader to \cite{Koc} for a detailed survey about the role of these operations and relations in the classification of two-dimensional topological quantum field theories. Let us note that a variant of Frobenius bialgebras called special Frobenius bialgebra is closely related to open-closed topological field theories \cite{LP} and conformal field theories \cite{FFRS}.

A classical example of Frobenius (bi)algebra comes from Poincaré duality. Let $M$ be an oriented connected closed manifold of dimension $n$. Let $[M]\in H_n(M;\mathbb{K})\cong H^0(M;\mathbb{K})\cong \mathbb{K}$ be the fundamental class of $[M]$. Then the cohomology ring $H^*(M;\mathbb{K})$ of $M$ inherits a structure of commutative and cocommutative Frobenius bialgebra of degree $n$ with the following data:

(i)the product is the cup product
\begin{eqnarray*}
\mu: H^kM\otimes H^lM & \rightarrow & H^{k+l}M \\
 x\otimes y & \mapsto & x\cup y
\end{eqnarray*}

(ii)the unit $\eta:\mathbb{K}\rightarrow H^0M\cong H_nM$ sends $1_{\mathbb{K}}$ on the fundamental class $[M]$;

(iii)the non-degenerate pairing is given by the Poincaré duality:
\begin{eqnarray*}
\beta: H^kM\otimes H^{n-k}M & \rightarrow & \mathbb{K}\\
 x\otimes y & \mapsto & <x\cup y,[M]>
\end{eqnarray*}
i.e the evaluation of the cup product on the fundamental class;

(iv) the coproduct $\Delta=(\mu\otimes id)\circ (id\otimes \gamma)$ where
\[
\gamma: \mathbb{K}\rightarrow \bigoplus_{k+l=n}H^kM\otimes H^lM
\]
is the dual copairing of $\beta$, which exists since $\beta$ is non-degenerate;

(v)the counit $\epsilon=<.,[M]>:H^nM\rightarrow \mathbb{K}$ i.e the evaluation on the fundamental class.
\end{example}

A natural question after looking at all these examples is the following: can we extract a common underlying pattern, analogue to representation theory of groups or to operad theory, which says that an algebraic structure of a given kind is an algebra over a corresponding combinatorial object ? A formalism that include algebras over operads as well as more general structures like Lie bialgebras and Frobenius bialgebras ?
We answer this question with the following definition, originally due to MacLane \cite{MLa}. A $\Sigma$-biobject is a double sequence $\{M(m,n)\in\Ch\}_{(m,n)\in\mathbb{N}^2}$ where each $M(m,n)$ is equipped with a right action of $\Sigma_{m}$ and a left action of $\Sigma_{n}$ commuting with each other.
\begin{defn}
A prop is a $\Sigma$-biobject endowed with associative horizontal composition products
\[
\circ_{h}:P(m_1,n_1)\otimes P(m_2,n_2)\rightarrow P(m_1+m_2,n_1+n_2),
\]
associative vertical composition products
\[
\circ_{v}:P(k,n)\otimes P(m,k)\rightarrow P(m,n)
\]
and maps $\mathbb{K}\rightarrow P(n,n)$ which are neutral for $\circ_v$ (representing the identity operations).
These products satisfy the exchange law
\[
(f_1\circ_h f_2)\circ_v(g_1\circ_h g_2) = (f_1\circ_v g_1)\circ_h(f_2\circ_v g_2)
\]
and are compatible with the actions of symmetric groups.

Morphisms of props are equivariant morphisms of collections compatible with the composition products.
\end{defn}
A fundamental example of prop is given by the following construction. To any complex $X$ we can associate an endomorphism prop $End_X$ defined by
\[
End_X(m,n)=Hom_{\Ch}(X^{\otimes m},X^{\otimes n}).
\]
The prop structure here is crystal clear: the actions of the symmetric groups are the permutations of the tensor powers, the vertical composition is the composition of homomorphisms and the horizontal one is the tensor product of homomorphisms.
\begin{defn}
A $P$-algebra on a complex $X$ is a prop morphism $P\rightarrow End_X$.
\end{defn}
That is, a $P$-algebra structure on $X$ is a collection of equivariant cochain morphisms
\[
\{P(m,n)\rightarrow Hom_{\Ch}(X^{\otimes m},X^{\otimes n})\}_{m,n\in\mathbb{N}}
\]
commuting with the vertical and horizontal composition products. Hence the formal operations of $P$ are sent to actual operations on $X$, and the prop structure of $P$ determines the relations satisfied by such operations.
\begin{rem}
MacLane's original definition is more compact: a prop $P$ in a closed symmetric monoidal category $\mathcal{C}$ as a symmetric monoidal category enriched in $\mathcal{C}$, with the natural integers as objects and the tensor product $\otimes$ defined by $m\otimes n=m+n$. A morphism of props is then an enriched symmetric monoidal functor. An algebra over a prop is an enriched symmetric monoidal functor $P\rightarrow\mathcal{C}$, and a morphism of algebras is an enriched symmetric monoidal transformation (see also \cite[Section 2.1]{Yal6} for the colored case).
\end{rem}

There is an adjunction between the category of $\Sigma$-biobjects and the category of props, with the right adjoint given by the forgetful functor and the left adjoint given by a free prop functor. Briefly, given a $\Sigma$-biobject $M$, the free prop $\mathcal{F}(M)$ on $M$ is defined by
\[
\mathcal{F}(M)(m,n)=\bigoplus_{G\in Gr(m,n)}(\bigotimes_{v\in Vert(G)}M(|In(v)|,|Out(v)|))_{Aut(G)}
\]
where
\begin{itemize}
\item The direct sums runs over the set $Gr(m,n)$ of directed graphs with $m$ inputs, $n$ outputs and no loops;

\item The tensor products are indexed by the sets $Vert(G)$ of vertices of such graphs $G$;

\item For each vertex $v$ of $G$, the numbers $|In(v)|$ and $|Out(v)|$ are respectively the number of inputs and the number of outputs of $v$;

\item These tensor products are mod out by the action of the group $Aut(G)$ of automorphisms of the graph $G$.
\end{itemize}
We refer the reader to \cite[Appendix A]{Frep} for more details about this construction.
Moreover, there is an obvious notion of ideal in a prop $P$, defined as a $\Sigma$-biobject $I$ such that $i\circ_v p\in I$ for $i\in I$ and $p\in P$, and $i\circ_h p\in I$ for $i\in I$ and $p\in P$. This means that each prop admits a \emph{presentation by generators and relations}, something particularly useful to describe an algebraic structure. For instance, all the operations $A^{\otimes n}\rightarrow A$ on an associative algebra $A$ induced by the algebra structure are entirely determined by a product $A\otimes A\rightarrow A$ and the associativity condition. Actually, the graphical presentations we gave in the examples above are exactly presentations of the corresponding props by generators and relations !
For instance, if we denote by $BiLie$ the prop of Lie bialgebra, we have
\[
BiLie=\mathcal{F}(M)/I
\]
where $M(2,1)=sgn_2\otimes\mathbb{K}.\xymatrix @R=0.5em@C=0.75em@M=0em{1\ar@{-}[dr] & & 2\ar@{-}[dl] \\ & \ar@{-}[d] & \\  & & }$, $M(1,2)=sgn_2\otimes\mathbb{K}.\xymatrix @R=0.5em@C=0.75em@M=0em{ & \ar@{-}[d] & \\ & \ar@{-}[dl] \ar@{-}[dr] & \\ 1 & & 2 }$ and $M(m,n)=0$ for $(m,n)\notin\{(2,1),(1,2)\}$ (recall here that $sgn_2$ is the signature representation of $\Sigma_2$). The ideal $I$ is generated by the graphs defining the relations in Example 4 (Jacobi, co-Jacobi, cocycle relation). A Lie bialgebra $g$ is then the datum of a prop morphism
\[
\{BiLie(m,n)\rightarrow Hom_{\Ch}(g^{\otimes m},g^{\otimes n})\}_{m,n\in\mathbb{N}}.
\]
According to the presentation of $BiLie$ by generators and relations, this prop morphism is completely determined by its values on the generators. That is, we send the generator $\xymatrix @R=0.5em@C=0.75em@M=0em{1\ar@{-}[dr] & & 2\ar@{-}[dl] \\ & \ar@{-}[d] & \\  & & }$ to a cochain map $[,]:g\otimes g\rightarrow g$, the generator $\xymatrix @R=0.5em@C=0.75em@M=0em{ & \ar@{-}[d] & \\ & \ar@{-}[dl] \ar@{-}[dr] & \\ 1 & & 2 }$ to a cochain map $\delta:g\rightarrow g\otimes g$, and the graphs of $I$ to zero. This implies that $[,]$ is a Lie bracket, $\delta$ a Lie cobracket and they satisfy moreover the cocycle relation.

Actually, for a wide range of algebraic structures, a well defined grafting operation on connected graphs is sufficient to parametrize the whole structure. Such a grafting is defined by restricting the vertical composition product of props to connected graphs. The unit for this connected composition product $\boxtimes_c$ is the $\Sigma$-biobject $I$ given by $I(1,1)=\mathbb{K}$ and $I(m,n)=0$ otherwise. The category of $\Sigma$-biobjects then forms a symmetric monoidal category $(Ch_{\mathbb{K}}^{\mathbb{S}},\boxtimes_c,I)$.
\begin{defn}
A dg properad $(P,\mu,\eta)$ is a monoid in $(Ch_{\mathbb{K}}^{\mathbb{S}},\boxtimes_c,I)$,
where $\mu$ denotes the product and $\eta$ the unit.
It is augmented if there exists a morphism of properads $\epsilon:P\rightarrow I$.
In this case, there is a canonical isomorphism $P\cong I\oplus\overline{P}$
where $\overline{P}=ker(\epsilon)$ is called the augmentation ideal of $P$.

Morphisms of properads are morphisms of monoids in $(Ch_{\mathbb{K}}^{\mathbb{S}},\boxtimes_c,I)$.
\end{defn}
Properads have also their dual notion, namely coproperads:
\begin{defn}
A dg coproperad $(C,\Delta,\epsilon)$ is a comonoid in $(Ch_{\mathbb{K}}^{\mathbb{S}},\boxtimes_c,I)$.
\end{defn}
As in the prop case, there exists a free properad functor $\mathcal{F}$ forming an adjunction
\[
\mathcal{F}:Ch_{\mathbb{K}}^{\mathbb{S}}\rightleftarrows Properad :U
\]
with the forgetful functor $U$.
There is an explicit construction of the free properad analogous to the free prop construction, but restricted to connected directed graphs instead of all directed graphs. Dually, there exists a cofree coproperad functor denoted $\mathcal{F}_c(-)$ having the same underlying $\Sigma$-biobject. There is also a notion of algebra over a properad similar to an algebra over a prop, since the endomorphism prop restricts to an endomorphism properad . Properads are general enough to encode a wide range of bialgebra structures such as associative and coassociative bialgebras, Lie bialgebras, Poisson bialgebras, Frobenius bialgebras for instance.
\begin{rem}
There is a free-forgetful adjunction between properads and props \cite{Val}.
\end{rem}

\section{Homotopy theory of (bi)algebras}

We already mentioned before the natural occurrence of ``relaxed'' algebraic structures, like $A_{\infty}$-algebras or $E_{\infty}$-algebras, in various situations where a given relation (associativity, commutativity) is satisfied only up to an infinite sequence of obstructions vanishing at the cohomology level. More generally, one can wonder how to set up a coherent framework to define what it means to ``relax'' a $P$-algebra structure, encompassing in particular the previous examples. Moreover, we will see later that deformation theory of differential graded $P$-algebras cannot be defined without working in the larger context of $P$-algebras up to homotopy (or homotopy $P$-algebras). This is due to the fact that the base category $\Ch$ itself manifests a non trivial homotopy theory. A natural way to define homotopy $P$-algebras is to resolve the prop $P$ itself by means of homotopical algebra. For this, we recall briefly that $\Ch$ has all the homotopical properties needed for our purposes, namely, it forms a cofibrantly generated symmetric monoidal model category. We refer the reader to to Hirschhorn \cite{Hir} and Hovey \cite{Hov} for a comprehensive treatment of homotopical algebra and monoidal model categories. The $\Sigma$-biobjects form a category of diagrams in $\Ch$ and inherit thus a cofibrantly generated model structure with pointwise weak equivalence and fibrations (the projective model structure). The free prop functor allows to transfer the projective model structure of $\Sigma$-biobjects along the free-forgetful adjunction:
\begin{thm}
(cf. \cite[Theorem 5.5]{Frep}) The category of dg props $Prop$ equipped with the
classes of componentwise weak equivalences and componentwise fibrations forms a cofibrantly generated model category.
\end{thm}
\begin{rem}
According to \cite{MV2}, the similar free-forgetful adjunction between $\Sigma$-biojects and dg properads equips dg properads with a cofibrantly generated model category structure with componentwise fibrations and weak equivalences.
\end{rem}
Hence we can define homotopy algebras over props as follows:
\begin{defn}
A homotopy $P$-algebra is a $P_{\infty}$-algebra, where $P_{\infty}\stackrel{\sim}{\rightarrow}P$ is a cofibrant resolution of $P$.
\end{defn}
Homotopy algebra structures appear naturally in plenty of topological and geometric situations, especially for transfer and realization problems:
\begin{itemize}
\item \textbf{Transfer problems}: given a quasi-isomorphism $X\stackrel{\sim}{\rightarrow}Y$, if $Y$ forms a $P$-algebra, then $X$ cannot inherits a $P$-algebra structure as well (since this is not a strict isomorphism) but rather a $P_{\infty}$-algebra structure. In the converse way, a choice of quasi-isomorphism $X\stackrel{\sim}{\rightarrow}H^*X$ from a complex to its cohomology allows to transfer any $P$-algebra structure on $X$ to a $P_{\infty}$-algebra structure on $H^*X$. That is, the data of a big complex with a strict structure can transferred to a smaller complex with a bigger structure up to homotopy.
\item \textbf{Realization problems}: A $P$-algebra structure on the cohomology $H^*X$ is induced by a finer $P_{\infty}$-algebra structure on $X$, which consists in a family of higher operations on cochains.
\end{itemize}
Let us name a few applications of such ideas:
\begin{itemize}
\item[$\bullet$] $A_{\infty}$-structures (associative up to homotopy) appeared very early in the study of loop spaces and monoidal categories (Stasheff's associahedra), and the $A_{\infty}$-structure induced on the singular cohomology of a topological space by the cochain-level cup product gives the higher Massey products. Such products are topological invariants, for instance the triple Massey product differentiate the borromean rings from the trivial link, even though their respective cohomologies are isomorphic as associative algebras.

\item[$\bullet$] $E_{\infty}$-structures (commutative up to homotopy) on ring spectra play a key role to encode cohomological operations in stable homotopy theory. Realization problems for such structures have been the subject of a consequent work by Goerss-Hopkins \cite{GoH}, following the idea of \cite{BDG} to study the homotopy type of the moduli space of all realizations on a given spectrum by decomposing it as the limit of a tower of fibrations, and determining the obstruction groups of the corresponding spectral sequence (which turns out to be André-Quillen cohomology groups).

\item[$\bullet$] The $E_{\infty}$-structure on singular cochains classifies the homotopy type of nilpotent spaces (see Sullivan over $\mathbb{Q}$, Mandell over $\mathbb{Z}$ and $\mathbb{F}_p$).

\item[$\bullet$] $L_{\infty}$-structures (Lie up to homotopy) encode the deformation theory of various algebraic, topological or geometric structures, a striking application being Kontsevich's deformation quantization of Poisson manifolds \cite{Ko2}.

\item[$\bullet$]In string topology, the homology of a loop space $\Omega M$ on a manifold $M$ is equipped with a natural Batalin-Vilkovisky algebra ($BV$-algebra) structure \cite{CS1}. On the other hand, the Hochschild cohomology of the singular cochains on $M$ is also a $BV$-algebra (extending the canonical Gerstenhaber algebra structure). In characteristic zero, when $M$ is a simply connected closed manifold, both are known to be isomorphic as $BV$-algebras \cite{FT}. It turns out that this structure lifts to a $BV_{\infty}$-structure on Hochschild cochains (a result called the cyclic Deligne conjecture). Homotopy $BV$-algebras are related not only to string topology but also to topological conformal field theories and vertex algebras \cite{GTV}.

\item[$\bullet$] Homotopy Gerstenhaber algebras, or equivalently $E_2$-algebras, are the natural structures appearing on Hochschild complexes by Deligne's conjecture, which has been generalized to the existence of $E_{n+1}$-algebra structures on higher Hochschild complexes of $E_n$-algebras. These results have applications to deformation quantization but also to factorization homology of manifolds and generalizations of string topology \cite{GTZ}. The proof of Deligne's conjecture relies on a transfer of structures combined with an obstruction theoretic method. Let us note that a bialgebra version of this conjecture obtained recently in \cite{GY} relies in particular on this ``transfer+obstruction'' method in the case of $E_3$-algebras and has applications to open problems in quantum group theory.
\end{itemize}
Moreover, homotopy algebra structures are the structures controlled by the cohomology theories of algebras, when one works in the dg setting. For instance, the Hochschild complex of a dg associative algebra $A$ controls (in a sense we will precise later) not the strict algebra deformations but the $A_{\infty}$ deformations of $A$.

However, there is a quite obvious problem in the definition of homotopy algebra we gave above. Indeed, it relies a priori on the choice of a resolution. For instance, two homotopy $P$-algebras could be weakly equivalent for a certain choice of $P_{\infty}$ but not for another choice. In order to make sense of this notion and of the various deformation theoretic, transfer and realization problems in which it naturally arises, we have to prove an invariance result for the homotopy theory of homotopy $P$-algebras:
\begin{thm}{\cite[Theorem 0.1]{Yal3}}
A weak equivalence $\varphi:P_{\infty}\stackrel{\sim}{\rightarrow}Q_{\infty}$ between cofibrant props induces an equivalence of $(\infty,1)$-categories
\[
\varphi^*:(Q_{\infty}-Alg,q-isos)\stackrel{\sim}{\rightarrow}(P_{\infty}-Alg,q-isos),
\]
where $(P_{\infty}-Alg,q-isos)$ is the $(\infty,1)$-category associated to the category of dg $P_{\infty}$-algebras with quasi-isomorphisms as weak equivalences.
\end{thm}
In the case of algebras over operads, this result is already known by using classical methods of homotopical algebra. A weak equivalence $\varphi:P\rightarrow Q$ of dg operads induces an adjunction
\[
\varphi_!:P_{\infty}-Alg\rightleftarrows Q_{\infty}-Alg:\varphi^*,
\]
where $\varphi^*$ is the functor induced by precomposition $P_{\infty}\rightarrow Q_{\infty}\rightarrow End_X$ and $\varphi_!$ is a certain coequalizer. The functor $\varphi^*$ is a right Quillen functor since weak equivalences and fibrations of algebras over operads are determined in complexes, so this is a Quillen adjunction. One can then prove that the unit and the counit of this adjunction are weak equivalences, hence the desired result (a Quillen equivalence induces an equivalence of the associated $(\infty,1)$-categories. We refer the reader to \cite[Chapter 16]{Fre3} for a detailed proof of this result. This method completely fails in the case of algebras over props for two reasons:
\begin{itemize}
\item Algebras over props are a priori not stable under all colimits, so the left adjoint $\varphi_!$ does not exist in general;

\item There is no free $P$-algebra functor, hence no way to transfer a model category structure from the one of cochain complexes (and by the previous point, the first axiom of model categories already fails).
\end{itemize}
To overcome these difficulties, one has to go through a completely new method based on the construction of a functorial path object of $P$-algebras and a corresponding equivalence of classification spaces proved in \cite{Yal2}, then an argument using the equivalences of several models of $(\infty,1)$-categories \cite{Yal3}. The equivalence of Theorem 2.4 is stated and proved in \cite{Yal3} as an equivalence of hammock localizations in the sense of Dwyer-Kan \cite{DK}.

Theorem 2.4 means that the notion of algebraic structure up to homotopy is coherent in a very general context, and in particular that transfer and realization problems make sense also for various kinds of bialgebras. Two motivating examples are the realizations of Poincaré duality of oriented closed manifolds as homotopy Frobenius algebra structures at the cochain level, and realizations of the Lie bialgebra structure on string homology at the chain level. Let us note that an explicit realization has been recently obtained in \cite{CFL} (with interesting relationships with symplectic field theory and Lagrangian Floer theory), using a notion of homotopy involutive Lie bialgebra which actually matches with the minimal model of the associated properad obtained in \cite{CMW} (see \cite[Remark 2.4]{CFL}). However, classification and deformation theory of such structures, as well as the potential new invariants that could follow, are still to be explored.

\section{Deformation theory and moduli problems in a derived framework}

\textbf{Geometric idea.} A common principle in algebraic topology and algebraic geometry is the following.
\begin{itemize}
\item In order to study a collection of objects (or structures) equipped with an equivalence relation, one construct a space (classifying space in topology, moduli space in geometry) whose points are given by this collection of objects and connected components are their equivalence classes.
\item The set of equivalences classes is not enough. Indeed, understanding the deformation theory of these objects amounts to study the infinitesimal deformations (formal neighbourhood) of the corresponding points on the moduli space. For this, one needs the existence of some tangent structure, thus the existence of a geometry on such a moduli space.
\item The deformation theory of a given point is then described by the associated formal moduli problem, which consists, roughly speaking, of a functor from augmented artinian cdgas to simplicial sets with nice gluing properties, so that its evaluation on an algebra $R$ is the space of $R$-deformations of this point.
\item One would like an algebraic description of this deformation theory in terms of deformation complexes and obstruction theory. For this, one has to move in the derived world and use Lurie's equivalence thm between formal moduli problems and dg Lie algebra. The corresponding dg Lie algebra is called the tangent Lie algebra.
\end{itemize}
In the two sections below, we describe some keys ideas to work out the construction above in a derived framework, and motivate the necessity to introduce these additional derived data.

\subsection{Derived algebraic geometry in a nutshell}

A usual geometric approach to moduli problems is to build an algebraic variety, scheme, or stack
parameterizing a given type of structures or objects (complex structures on a Riemann surface, vector bundles of fixed rank...). However, the usual stacks theory shows its limits when one wants to study families of objects related by an equivalence notion weaker than isomorphisms (for instance, complexes of vector bundles) and capture their full deformation theory on the tangent spaces. Derived algebraic geometry is a conceptual framework to solve such problems, that can be seen as a homotopical perturbation or thickening of algebraic geometry \cite{TV}.

Recall that as a ringed space, a usual scheme is a couple $(X,\mathcal{O}_X$, where $X$ is a topological space and $\mathcal{O}_X$ a sheaf of commutative algebras over $X$ called the structural sheaf of the scheme. That is, \emph{schemes are structured spaces locally modelled by commutative algebras}. From the ``functor of points'' perspective, schemes are sheaves $Aff\rightarrow Set$ on the category $Aff$ of affine schemes, which is the opposite category of the category $Com_{\mathbb{K}}$ of commutative algebras: they are functors $Com_{\mathbb{K}}\rightarrow Set$ satisfying a gluing condition (also called descent condition) with respect to a specified collection of families of maps in $Com_{\mathbb{K}}$ called a Grothendieck topology on $Com_{\mathbb{K}}$. The notion of Grothendieck topology can be seen as a categorical analogue of the notion of covering of a topological space, and like sheaves on a topological space, we want sheaves on a given category to satisfy a gluing condition along the ``coverings'' given by this Grothendieck topology. Stack theory goes one step further, replacing $Set$ by the $2$-category of groupoids $Grpd$. Stacks are then functors $Com_{\mathbb{K}}\rightarrow Grpd$ satisfying a $2$-categorical descent condition (gluing on objects of the groupoids and compatible gluing on sets of isomorphisms between these objects). A motivation for such a complicated generalization of scheme theory is to handle all the interesting moduli problems that cannot be represented by a moduli space in the category of schemes, due to the fact that the families of objects parametrized by this moduli problem have \emph{non trivial automorphisms} (consider for instance fiber bundles on a variety).

To give a geometric meaning and good properties for such moduli spaces, one has to go further and work with geometric stacks, a subcategory of stacks which can be obtained by gluing (taking quotients of) representable stacks along a specified class $\textbf{P}$ of maps. An important example of Grothendieck topology is the étale topology. In this topology, the geometric stacks obtained by choosing for $\textbf{P}$ the class of étale maps are the Deligne-Mumford stacks, and the geometric stacks obtained by choosing for $\textbf{P}$ the class of smooth maps are the Artin stacks. To satisfy the corresponding conditions forces the points of such a stack to have ``not too wild'' automorphism groups: the points of a Deligne-Mumford stack have finite groups of automorphisms (the historical example motivating the introduction of stack theory is the moduli stack of stable algebraic curves), and Artin stacks allow more generally algebraic groups of automorphisms (for example a quotient of a scheme by the action of an algebraic group).

A derived scheme is a couple $S=(X,\mathcal{O}_X$, where $X$ is a topological space and $\mathcal{O}_X$ a sheaf of differential graded commutative algebras over $X$, such that $t_0S=(X,H^0\mathcal{O}_X)$ (the zero truncation of $S$) is a usual scheme and the $H^{-i}\mathcal{O}_X$ are quasi-coherent modules over $H^0\mathcal{O}_X$. That is, \emph{derived schemes are structured spaces locally modelled by cdgas}. Using the ``functor of points'' approach, we can present derived geometric objects in the diagram
\[
\xymatrix{
Com_{\mathbb{K}}\ar[r]\ar[dd]\ar[dr]\ar[ddr] & Set \ar[d] \\
 & Grpd \ar[d] \\
CDGA_{\mathbb{K}} \ar[r] & sSet
}
\]
where $CDGA_{\mathbb{K}}$ is the $\infty$-category of non-positively graded commutative differential graded algebras, $Set$ the category of sets, $Grpds$ the ($2$-)category of groupoids and $sSet$ the $\infty$-category of simplicial sets ($\infty$-groupoids).
\begin{itemize}
\item Schemes are sheaves $Com_{\mathbb{K}}\rightarrow Set$ over the category of affine schemes (the opposite category of $Com_{\mathbb{K}}$) for a choice of Grothendieck topology.
\item Stacks are ``sheaves'' $Com_{\mathbb{K}}\rightarrow Grpd$ for a $2$-categorical descent condition, and landing in groupoids allows to represent moduli problems for which objects have
\textbf{non-trivial automorphisms}.
\item Higher stacks are ``sheaves up to homotopy'' $Com_{\mathbb{K}}\rightarrow sSet$, and landing in simplicial sets allows to represent moduli problems for which objects are related by
\textbf{weak equivalences} instead of isomorphisms.
\item Derived stacks are ``sheaves up to homotopy'' $CDGA_{\mathbb{K}}\rightarrow sSet$ over the $\infty$-category of non-positively graded cdgas (in the cohomological convention) with a choice of Grothendieck topology on the associated homotopy category. They capture the \textbf{derived data} (obstruction theory via (co)tangent complexes, non-transverse intersections, $K$-theoretic virtual fundamental classes \cite[Section 3]{Toe}) and convey richer geometric structures (shifted symplectic structures for instance \cite{PTVV}).
\end{itemize}
It is important to precise that, to get derived stacks with geometric properties, we have to restrict to a sub-$\infty$-category of these, called derived Artin stacks. Derived $1$-Artin stacks are geometric realizations of smooth groupoids objects in derived affine schemes, and derived $n$-Artin stack are recursively defined as the geometric realization of smooth groupoid object in derived $n-1$-Artin stacks. An alternate way is to define $n$-Artin stacks as smooth $n$-hypergroupoid objects in derived affine schemes \cite{Pri2}. This is the natural generalization, in the derived setting, of the geometric stacks we mentionned earlier: we obtain them by gluing representables along smooth maps, and this gluing is defined as the realization of a (higher) ``groupoid-like'' object. Such stacks are also said to be $n$-geometric. Derived Artin stacks admit cotangent complexes, an associated obstruction theory and various properties for which we refer the reader to \cite{Toe}. Concerning in particular the obstruction theory, the cotangent complex of a derived $n$-Artin stack is cohomologically concentrated in degrees $]-\infty;n]$. If the derived Artin stack $X$ is locally of finite presentation, then it admits a tangent complex (the dual of the cotangent complex in the $\infty$-category $L_{qcoh}(X)$ of quasi-coherent complexes over $X$) cohomologically concentrated in degrees $[-n;\infty[$. The geometric meaning of the cohomological degree is the following: at a given point $x$ of $X$, the cohomology of the tangent complex in positive degrees controls the \emph{obstruction theory} of $x$ (extensions of infinitesimal deformations to higher order deformations), the $0^{th}$-cohomology group is the group of \emph{equivalence classes of infinitesimal deformations}, and the cohomology of the tangent complex in negative degrees controls the \emph{(higher) symmetries} of $x$ (the homotopy type of its automorphisms is bounded by $n$). This last part generalizes to derived geometry the idea of the usual theory of algebraic stacks, that we have to control the automorphisms of the points to get a nice geometric object.
\begin{rem}
Derived Artin stacks satisfy the ``geometricity'' condition for a derived analogue of the class of smooth maps. Similarly, one can define derived Deligne-Mumford stacks by a geometricity condition for a derived analogue of the class of étale maps.
\end{rem}

To illustrate this homotopical enhancement of algebraic geometry, let us give some interesting examples.
\begin{example}
Let $X$ and $Y$ be two subvarieties of a smooth variety $V$. Their intersection is said to be transverse if and only if for every point $p\in X\cap Y$, we have $T_pV=T_pX+T_pY$ where $T_p$ is the tangent space at $p$. This means that $X\cap Y$ is still a subvariety of $V$. Transverse intersections are very useful:
\begin{itemize}
\item[$\bullet$] In algebraic topology, to define the intersection product $[X].[Y]=[X\cap Y]$ on the homology $H_*M$ of a manifold $M$ (classes being represented by submanifolds $X,Y$ of $M$).

\item[$\bullet$] In algebraic geometry, classes represented by subvarieties are called algebraic classes, and the formula of the intersection product above equip algebraic classes with a ring structure. This is called the Chow ring.
\end{itemize}
It is thus natural to ask what happens when intersections are not transverse. The idea is to deform $X$ to another subvariety $X'$ and $Y$ to another subvariety $Y'$ such that $X'$ and $Y'$ intersect transversely, and to define $[X].[Y]=[X'\cap Y']$. The drawback is that $X\cap Y$ is not a geometric object anymore but just a homology class.

Another natural question is to count multiplicity (in some sense, the ``degree of tangency'') of non transverse intersections. For example, consider $X=\{y=0\}$ a line tangent to $Y=\{y-x^2=0\}$ the parabola in the affine plane, and look at the intersection point $p=(0,0)$ of $X$ and $Y$. If we deform this situation to a generic case by moving the line along the parabola, the line intersects the parabola at two distinct points. This means that the multiplicity of $p$ is $2$. In general, the multiplicity of the intersection of two subvarieties $X$ and $Y$ at a generic point $p$ is given by Serre's intersection formula
\begin{eqnarray*}
I(p;X,Y) & = & \sum_i(-1)^i dim_{\mathcal{O}_{V,p}}(Tor^{\mathcal{O}_{V,p}}_i(\mathcal{O}_{X,p},\mathcal{O}_{Y,p})) \\
 & = & dim(\mathcal{O}_{X,p}\otimes_{\mathcal{O}_{V,p}}\mathcal{O}_{Y,p})+\textrm{correction terms}
\end{eqnarray*}
where $\mathcal{O}_{V,p}$ is the stalk of $\mathcal{O}_V$ at $p$, and $\mathcal{O}_{X,p},\mathcal{O}_{Y,p}$ are $\mathcal{O}_{V,p}$-modules for the structures induced by the inclusions $X\hookrightarrow V,Y\hookrightarrow V$. In certain cases, the multiplicity is determined by the dimension of $\mathcal{O}_{X,p}\otimes_{\mathcal{O}_{V,p}}\mathcal{O}_{Y,p}$, but in general this is not sufficient and we have to introduce correction terms given by the derived functors $Tor$ with no geometric meaning.

Non transverse intersections have a natural geometric construction in derived geometry. The idea is to realize $X\cap Y$ as a derived scheme by using a derived fiber product
\[
X\times^h_V Y=(X\cap Y, \mathcal{O}_{X\times^h_V Y}=\mathcal{O}_{X}\otimes_{\mathcal{O}_{V}}^{\mathbb{L}}\mathcal{O}_{Y})
\]
where $\otimes^{\mathbb{L}}$ is the left derived tensor product of sheaves of cdgas and $\otimes_{\mathcal{O}_{V}}^{\mathbb{L}}$ is the left derived tensor product of dg $\mathcal{O}_V$-modules. Then
\[
I(p;X,Y) = \sum_i(-1)^i dim(H^{-i}\mathcal{O}_{X\times^h_V Y}),
\]
that is, the intersection number naturally and geometrically arises as the Euler characteristic of the structure sheaf of the derived intersection. In a sentence, \emph{the transversality failure is measured by the derived part of the structure sheaf}.
\end{example}
\begin{example}
Another kind of application is Kontsevich's approach to Gromov-Witten theory in symplectic topology and algebraic geometry (which has also applications in string theory). On the algebraic geometry side, the problem is the following. When we want to count the intersection points of two curves in $\mathbb{P^2}$, we use intersection theory on $\mathbb{P}^2$ and Bezout theorem. More generally, one could wonder how to count rational curves of a given degree in $\mathbb{P}^N$ that intersect a given number of points $p_1,\cdots,p_n$, or replace $\mathbb{P}^N$ by a more general variety $X$. The idea to address this question is to define a moduli space of such curves and do intersection theory on this moduli space. But for this, one has to define a moduli space with good geometric properties, a constraint that leads to the notion of stable map. Let $C$ be a curve of genus $g$ and degree $d$ with marked points $p_1,\cdots,p_n$. A stable map is a map $f:C\rightarrow X$ satisfying an additional ``stability condition'' we do not precise here. Counting rational curves of genus $g$ and degree $d$ in $X$ passing through $n$ fixed points $x_1,\cdots,x_n$ of $X$ amounts to count such stable maps, and this defines the Gromov-Witten invariants of $X$. A classical idea is to define an invariant by integrating some function on the appropriate moduli space (via intersection theory). Here, this is the moduli space of stable maps $\overline{\mathcal{M}_{g,n}}(X,d)$. In the case $X=\mathbb{P}^n$, this is a smooth and compact Deligne-Mumford stack. In the general case of a smooth proper variety, the moduli space $\overline{\mathcal{M}_{g,n}}(X,d)$ is not smooth anymore and this is a major trouble.

Indeed, we would like to define Gromov-Witten invariants by
\begin{eqnarray*}
GW_d(x_1,\cdots,x_n) & = & \int_{\overline{\mathcal{M}_{g,n}}(X,d)}ev_1^*[x_1]\cup\cdots\cup ev_n^*[x_n] \\
 & =& <[\overline{\mathcal{M}_{g,n}}(X,d)],ev_1^*[x_1]\cup\cdots\cup ev_n^*[x_n]>
\end{eqnarray*}
where $ev_i:\overline{\mathcal{M}_{g,n}}(X,d)\rightarrow X,f\mapsto f(p_i)$ is the evaluation map at the $i^{th}$ marked point of curves, the class $[x_i]$ is the cohomology class associated to the homology class of the point $x_i$ by Poincaré duality, and $<,>$ is the Poincaré duality pairing. Intuitively, the class $ev_1^*[x_i]$ represents curves in $X$ whose $i^{th}$ marked point coincide (up to deformation of the curve) with $x_i$, that is, equivalences classes of curves passing through $x_i$. The product $ev_1^*[x_1]\cup\cdots\cup ev_n^*[x_n]$ then correspond to equivalence classes of curves passing through $x_1,\cdots,x_n$, and counting such curves amounts to pair it along the fundamental class $[\overline{\mathcal{M}_{g,n}}(X,d)]$ of $\overline{\mathcal{M}_{g,n}}(X,d)$. And this is the problem: there is no such thing as a ```fundamental class of $\overline{\mathcal{M}_{g,n}}(X,d)$'', since $\overline{\mathcal{M}_{g,n}}(X,d)$ is not smooth.

Briefly, Kontsevich's idea is to see $\overline{\mathcal{M}_{g,n}}(X,d)$ as a ``derived space'' (i.e. equipped with a differential graded structure sheaf), that is, to make $\overline{\mathcal{M}_{g,n}}(X,d)$ formally behave like a smooth space by replacing the tangent spaces by tangent complexes.
Then one associates to its dg sheaf a ``virtual fundamental class'' $[\overline{\mathcal{M}_{g,n}}(X,d)]^{vir}$, generalizing the fundamental class of smooth objects to singular objects (by taking the Euler characteristic of this dg sheaf in $K$-theory, and sending this $K$-theory virtual class to a class in the Chow ring of $\overline{\mathcal{M}_{g,n}}(X,d)$, thanks to the existence of a Chern character). This allows to properly define
\[
GW_d(x_1,\cdots,x_n) = <[\overline{\mathcal{M}_{g,n}}(X,d)]^{vir},ev_1^*[x_1]\cup\cdots\cup ev_n^*[x_n]>.
\]
\end{example}
\begin{example}
Another very interesting application is the possibility to define a derived version of character varieties. Let $M$ be a smooth manifold and $G$ a Lie group (or an algebraic group). We know that a $G$-local system on $M$ is given by a $G$-bundle with flat connection, and those bundles are equivalent to representations $\pi_1M\rightarrow G$ by the Riemann-Hilbert correspondence. The variety of $G$-characters of $M$ is defined by
\[
Loc_G(M)=Hom(\pi_1M,G)/G
\]
where $G$ acts by conjugation. This is the moduli space of $G$-local systems on $M$. Character varieties are of crucial importance in various topics of geometry and topology, including
\begin{itemize}
\item[$\bullet$] Teichmüller geometry: for a Riemann surface $S$, the variety $Loc_{SL_2}(S)$ contains the Teichmüller space of $S$ as a connected component.

\item[$\bullet$] Low dimensional topology: for $dim(M)=3$, the variety $Loc_G(M)$ is related to quantum Chern-Simons invariants of $M$ (their are various conjectures about how the properties of $Loc_G(M)$ could determine the behaviour of the $3$-TQFT associated to $G$ and $M$ and associated invariant such as the colored Jones polynomial).
\end{itemize}
However, this is generally a highly singular object, and one would like to apply the principle shown in the previous example: treat this singular object as a smooth object in a derived framework. To formalize this idea, one defines a derived stack
\[
RLoc_G(M)=Map(Betti(M),BG)
\]
where $Betti(M)$ is the Betti stack of $M$, $BG$ is the derived classifying stack of $M$ and $Map$ is the internal mapping space in the $\infty$-category of derived stacks \cite{Toe}. This new object satisfies the following important properties:
\begin{itemize}
\item[$\bullet$] Its zero truncation gives the usual character variety
\[
T_0RLoc_G(M)=Loc_G(M).
\]

\item[$\bullet$] The tangent complex over a point computes the cohomology of $M$ with coefficients in the associated $G$-local system.

\item[$\bullet$] There is a nice new geometric structure appearing on such objects, which is typically of derived nature: it possesses a canonical $2-dim(M)$-shifted symplectic structure \cite{PTVV}. Briefly, shifted symplectic structures are the appropriate generalization of symplectic structure from smooth manifolds to derived stacks. Here, since tangent spaces are complexes, differential forms come with a cohomological degree in addition to their weight. An $n$-shifted symplectic structure is a cohomology class of degree $n$ in the de Rham complex of closed $2$-forms satisfying a weak non-degeneracy condition: for every point $Spec(\mathbb{K})\rightarrow X$, it induces a quasi-isomorphism $\mathbb{T}_{X/\mathbb{K}}\stackrel{\sim}{\rightarrow}\mathbb{L}_{X/\mathbb{K}}[-n]$ between the tangent complex and the shifted cotangent complex.
\end{itemize}
If $X$ is a smooth manifold, a $0$-shifted symplectic structure on $X$ is a usual symplectic structure. Let us note that in the case of a surface, the $0$-shifted symplectic form $RLoc_G(M)$ restricts to Goldman's symplectic form on the smooth locus of $Loc_G(M)$ \cite{Gol}, so this is really an extension of Goldman's form to the whole moduli space.
\end{example}

Finally, to come back to the main topic of our survey and to motivate a bit the use of homotopy theory for moduli spaces of algebraic structures, let us see on a very simple example what happens if we build such a space with usual algebraic geometry:
\begin{example}
Let $V$ be a vector space of dimension $n$, and let us consider a basis $\{e_1,\cdot,e_n\}$ of $V$. An associative product on $V$ is a linear map $\mu:V\otimes V\rightarrow V$ satisfying the associativity condition, hence it is determined by its values on the basis vectors
\[
\mu (e_i,e_j)=\sum_{k=1}^nc_{ij}^ke_k,
\]
where the $c_{ij}^k$'s satisfy moreover a certain set of relations $R$ determined by the associativity of $\mu$.
We can build an affine scheme whose $\mathbb{K}$-points are the associative algebra structures on $V$: its ring of functions is simply given by  $A=\mathbb{K}[c_{ij}^k]/(R)$. But we would like to classify such structures up to isomorphism, hence up to base change in $V$. For this, we have to mod out by the action of $GL_n$ on $V$. In order to have a well defined quotient of $Spec(A)$ by $GL_n$, we take the quotient stack $[Spec(A)/GL_n]$ as our moduli space of associative algebra structures up to isomorphisms.

Now let $R$ be an associative algebra with underlying vector space $V$, which represents a $\mathbb{K}$-point of $V$ (given by the orbit of the action of $GL(V)$ on $R$). Then the truncated tangent complex $\mathbb{T}_R$ of $[Spec(A)/GL_n]$ over the orbit of $R$ is given by a map
\[
d\psi : gl(V)\rightarrow T_RSpec(A),
\]
where $gl(V)$ is the Lie algebra of $GL(V)$ (the Lie algebra of matrices with coefficients in $V$) sitting in degree $-1$, and $T_RSpec(A)$ is the tangent space of $Spec(A)$ over $R$, sitting in degree $0$. This map is the tangent map of the scheme morphism
\[
\phi : GL(V)\rightarrow Spec(A)
\]
which sends any $f\in GL(V)$ to $f.R$, the action of $f$ on $R$, defined by transferring the algebra structure of $R$ along $f$. This is what one should expect for the tangent complex: two associative algebra structures are equivalent if and only if they are related by the action of $GL(V)$ (also called action of the ``gauge group'').
We then get $H^{-1}\mathbb{T}_R=End_{alg}(R)$ (the Lie algebra of algebra endomorphisms of $R$, tangent to $Aut_{alg}(R)$) and $H^0\mathbb{T}_R=HH^2(R,R)$ the second Hochschild cohomology group of $R$. Let us note that this computation is a very particular case of \cite[Theorem 5.6]{Yal4}. The group $HH^2(R,R)$ classifies equivalence classes of infinitesimal deformations of $R$. In particular, if $HH^2(R,R)=0$ then the algebra $R$ is rigid, in the sense that any infinitesimal deformation of $R$ is equivalent to the trivial one.
\end{example}
The construction above has two main drawbacks. First, the tangent complex does not give us any information about the obstruction theory of $R$, for instance, obstruction groups for the extension of infinitesimal deformations to formal ones. Second, in the differential graded case this construction does not make sense any more, and gives no way to classify structures up to quasi-isomorphisms.

\subsection{Derived formal moduli problems}

Formal moduli problems arise when one wants to study the infinitesimal deformation theory of a point $x$ of a given moduli space $X$ (variety, scheme, stack, derived stack) in a formal neighbourhood of this point (that is, the formal completion of the moduli space at this point). Deformations are parametrized by augmented artinian rings, for example $\mathbb{K}[t]/(t^2)$ for infinitesimal deformations of order one, or $\mathbb{K}[t]/(t^n)$ for polynomial deformations of order $n$. The idea is to pack all the possible deformations of $x$ in he datum of a deformation functor
\[
Def_{X,x}:Art_{\mathbb{K}}^{aug}\rightarrow Set
\]
from augmented artinian algebras to sets, sending an artinian algebra $R$ to the set of equivalence classes of $R$-deformations of $x$, that is, equivalence classes of lifts
\[
\xymatrix{
Spec(R)\ar[r] & X \\
Spec(\mathbb{K})\ar[u] \ar[ur]_-x & 
}
\]
(where the morphism $Spec(\mathbb{K})\rightarrow Spec(R)$ is induced by the augmentation $R\rightarrow\mathbb{K}$).
These are nothing but the fiber of the map $X(R)\rightarrow X(\mathbb{K})$ induced by the augmentation $R\rightarrow\mathbb{K}$ and taken over the base point $x$. Later on, several people realized that one could use Lie theory of dg Lie algebra to describe these deformation functors. Precisely, given a dg Lie algebra $g$, we consider the functor
\begin{eqnarray*}
Def_g:Art_{\mathbb{K}}^{aug} & \rightarrow & Set \\
R & \longmapsto & MC(g\otimes_{\mathbb{K}} m_R)
\end{eqnarray*}
where $m_R$ is the maximal ideal of $g$ and  $MC(g\otimes_{\mathbb{K}} m_R)$ is the set of Maurer-Cartan elements of the dg Lie algebra $g\otimes_{\mathbb{K}} m_R$, that is, elements $x$ of degree $1$ satisfying the Maurer-Cartan equation $dx+\frac{1}{2}[x,x]=0$. The functor $Def_g$ is a formal moduli problem called the deformation functor or deformation problem associated to $g$. This characterization of formal moduli problems arose from unpublished work of Deligne, Drinfed and Feigin, and was developed further by Goldman-Millson, Hinich, Kontsevich, Manetti among others. Defining deformation functors via dg Lie algebras led to striking advances, for instance in the study of representations of fundamental groups of varieties \cite{GM,Sim} and in deformation quantization of Poisson manifolds \cite{Ko2}.

It turned out that all known deformation problems related to moduli spaces in geometry were of this form, which led these people to conjecture that there should be a general correspondence between formal moduli problems and dg Lie algebras. However, there was no systematic recipe to build a dg Lie algebra from a given moduli problem (the construction above is the converse direction of this hypothetical equivalence), and even worse, different dg Lie algebras could represent the same moduli problem. Moreover, the obstruction theory associated to a moduli problem, given by the positive cohomology groups of its Lie algebra, has no natural interpretation in terms of the deformation functor. Indeed, deformation theory is actually of derived nature. For instance, if we want to study the extension of polynomial deformations of order $n$ to order $n+1$, we have to study the properties of the natural projection $\mathbb{K}[t]/(t^{n+1})\rightarrow \mathbb{K}[t]/(t^n)$ and under which conditions the induced map $X(\mathbb{K}[t]/(t^{n+1}))\rightarrow X(\mathbb{K}[t]/(t^n))$ is surjective, or bijective. This projection actually fits in a \emph{homotopy pullback} (not a strict pullback) of \emph{augmented dg artinian algebras} (not augmented commutative algebras in $\mathbb{K}$-modules)
\[
\xymatrix{
\mathbb{K}[t]/(t^{n+1})\ar[r]\ar[d] & \mathbb{K}[t]/(t^n)\ar[d]\\
\mathbb{K}\ar[r] & \mathbb{K}[\epsilon]/(\epsilon^2)
},
\]
where $\epsilon$ is of cohomological degree $1$ (not $0$).
If we could define formal moduli problems in this dg setting, we would like to apply the formal moduli problem $X_x$, associated to a given point $x$ of a moduli space $X$, to the diagram above to get a fiber sequence
\[
X_x(\mathbb{K}[t]/(t^{n+1}))\rightarrow X_x(\mathbb{K}[t]/(t^n))\rightarrow X(\mathbb{K}[\epsilon]/(\epsilon^2))
\]
and study the obstruction theory by understanding $X(\mathbb{K}[\epsilon]/(\epsilon^2))$ in an algebraic way.

These problems hints towards the necessity to introduce some homotopy theory in the study of formal moduli problems. For this, one replaces augmented artinian algebras $Art_{\mathbb{K}}^{aug}$ by augmented dg artinian algebras $dgArt_{\mathbb{K}}^{aug}$, and sets $Set$ by simplicial sets $sSet$:
\begin{defn}
A derived formal moduli problem is a functor $F:dgArt_{\mathbb{K}}^{aug}\rightarrow sSet$ from augmented artinian commutative differential graded algebras to simplicial sets, such that

1. We have an equivalence $F(\mathbb{K})\simeq pt$.

2. The functor $F$ sends quasi-isomorphisms of cdgas to weak equivalences of simplicial sets.

3. Let us consider a homotopy pullback of augmented dg artinian algebras
\[
\xymatrix{
A\ar[r]\ar[d] & B\ar[d] \\
C\ar[r] & D }
\]
and suppose that the induced maps $H^0C\rightarrow H^0D$ and $H^0B\rightarrow H^0D$ are surjective. Then $F$ sends this homotopy pullback to a homotopy pullback of simplicial sets.
\end{defn}
Formal moduli problems form a full sub-$\infty$-category noted $FMP_{\mathbb{K}}$ of the $\infty$-category of simplicial presheaves over augmented artinian cdgas. To explicit the link with derived algebraic geometry, the formal neighbourhood of a point $x$ in a derived stack $X$ (formal completion of $X$ at $x$) gives the derived formal moduli problem $X_x$ controlling the deformation theory of $x$. Given an artinian algebra $R$ with augmentation $\epsilon:R\rightarrow\mathbb{K}$, the homotopy fiber
\[
X_x(R)=hofib(X(\epsilon):X(R)\rightarrow X(\mathbb{K}))
\]
taken over the $\mathbb{K}$-point $x$ is the space of $R$-deformations of $X$, and equivalence classes of $R$-deformations are determined by $\pi_0X_x(R)$. In particular, applying $X_x$ to the homotopy pullback
\[
\xymatrix{
\mathbb{K}[t]/(t^{n+1})\ar[r]\ar[d] & \mathbb{K}[t]/(t^n)\ar[d]\\
\mathbb{K}\ar[r] & \mathbb{K}[\epsilon]/(\epsilon^2)
},
\]
we get a homotopy fiber sequence of spaces
\[
X_x(\mathbb{K}[t]/(t^{n+1}))\rightarrow X_x(\mathbb{K}[t]/(t^n))\rightarrow X(\mathbb{K}[\epsilon]/(\epsilon^2)),
\]
hence a fiber sequence
\[
\pi_0X_x(\mathbb{K}[t]/(t^{n+1}))\rightarrow \pi_0X_x(\mathbb{K}[t]/(t^n))\rightarrow \pi_0X_x(\mathbb{K}[\epsilon]/(\epsilon^2))\cong H^1\mathfrak{g}_{X_x},
\]
where $\mathfrak{g}_{X_x}$ is the tangent Lie algebra of the formal moduli problem $X_x$. We can take equivalently the cohomology of the shifted tangent complex $\mathbb{T}_{X,x}[-1]$ of the stack $X$ at $x$.
\begin{rem}
Actually, as proved in \cite{Hen}, for any derived Artin stack $X$ locally of finite presentation (so that we can dualize the cotangent complex to define the tangent complex), there exists a quasi-coherent sheaf $\mathfrak{g}_X$ of $\mathcal{O}_X$-linear dg Lie algebras over $X$ such that
\[
\mathfrak{g}_X\simeq \mathbb{T}_{X/\mathbb{K}}[-1]
\]
in the $\infty$-category $L_{qcoh}(X)$ of quasi-coherent complexes over $X$, where $\mathbb{T}_{X/\mathbb{K}}$ is the global tangent complex of $X$ over $\mathbb{K}$. Pulling back this equivalence along a point $x:Spec(\mathbb{K})\rightarrow X$, we get a quasi-isomorphism $\mathfrak{g}_{X_x}\simeq \mathbb{T}_{X,x}[-1]$. The sheaf $\mathfrak{g}_X$ thus encodes the family of derived formal moduli problems parametrized by $X$ which associates to any point of $X$ its deformation problem (the formal completion of $X$ at this point).
\end{rem}
The rigorous statement of an equivalence between derived formal moduli problems and dg Lie algebras was proved independently by Lurie in \cite{Lur0} and by Pridham in \cite{Pri}:
\begin{thm}[Lurie, Pridham]
The $\infty$-category $FMP_{\mathbb{K}}$ of derived formal moduli problems over $\mathbb{K}$ is equivalent to the $\infty$-category $dgLie_{\mathbb{K}}$ of dg Lie $\mathbb{K}$-algebras.
\end{thm}
Moreover, one side of the equivalence is made explicit, and is equivalent to the nerve construction of dg Lie algebras studied thoroughly by Hinich in \cite{Hin0}. The homotopy invariance of the nerve relies on nilpotence conditions on the dg Lie algebra. In the case of formal moduli problems, this nilpotence condition is always satisfied because one tensors the Lie algebra with the maximal ideal of an augmented artinian cdga. In this article , what we will call moduli problems are actually derived moduli problems.

\subsubsection{Extension to $L_{\infty}$-algebras}

Certain deformation complexes of interest are not strict Lie algebras but homotopy Lie algebras, that is $L_{\infty}$-algebras. Their is a strictification theorem for homotopy Lie algebras (more generally, for dg algebras over any operad when $\mathbb{K}$ is of characteristic zero), so any $L_{\infty}$-algebra is equivalent to a dg Lie algebra, but this simplification of the algebraic structure goes with an increased size of the underlying complex, which can be very difficult to explicit. This is why one would like the theory of derived formal moduli problems to extend to $L_{\infty}$-algebras, and fortunately it does. There are two equivalent definitions of an $L_{\infty}$-algebra:
\begin{defn}
(1) An $L_{\infty}$-algebra is a graded vector space $g=\{g_n\}_{n\in\mathbb{Z}}$ equipped with maps
$l_k:g^{\otimes k}\rightarrow g$ of degree $2-k$, for $k\geq 1$, satisfying the following properties:
\begin{itemize}
\item $l_k(...,x_i,x_{i+1},...)=-(-1)^{|x_i||x_{i+1}|}l_k(...,x_{i+1},x_i,...)$
\item for every $k\geq 1$, the generalized Jacobi identities
\[
\sum_{i=1}^k\sum_{\sigma\in Sh(i,k-i)}(-1)^{\epsilon(i)}l_k(l_i(x_{\sigma(1)},...,x_{\sigma(i)}),x_{\sigma(i+1)},...,x_{\sigma(k)})=0
\]
where $\sigma$ ranges over the $(i,k-i)$-shuffles and
\[
\epsilon(i) = i+\sum_{j_1<j_2,\sigma(j_1)>\sigma(j_2)}(|x_{j_1}||x_{j_2}|+1).
\]
\end{itemize}

(2) An $L_{\infty}$-algebra structure on a graded vector space $g=\{g_n\}_{n\in\mathbb{Z}}$ is a
coderivation $Q:\hat{Sym}^{\bullet\geq 1}(g[1])\rightarrow \hat{Sym}^{\bullet\geq 1}(g[1])$ of degree $1$ of the cofree cocommutative coalgebra $\hat{Sym}^{\bullet\geq 1}(g[1])$ such that $Q^2=0$.
\end{defn}
The bracket $l_1$ is actually the differential of $g$ as a cochain complex. When the brackets $l_k$ vanish
for $k\geq 3$, then one gets a dg Lie algebra.
The dg algebra $C^*(g)$ obtained by dualizing the dg coalgebra of (2) is called the Chevalley-Eilenberg algebra of $g$.

A $L_{\infty}$ algebra $g$ is filtered if it admits a decreasing filtration
\[
g=F_1g\supseteq F_2g\supseteq...\supseteq F_rg\supseteq ...
\]
compatible with the brackets: for every $k\geq 1$,
\[
l_k(F_rg,g,...,g)\in F_rg.
\]
We suppose moreover that for every $r$, there exists an integer $N(r)$ such that $l_k(g,...,g)\subseteq F_rg$
for every $k>N(r)$.
A filtered $L_{\infty}$ algebra $g$ is complete if the canonical map $g\rightarrow lim_rg/F_rg$ is an isomorphism.

The completeness of a $L_{\infty}$ algebra allows to define properly the notion of Maurer-Cartan element:
\begin{defn}
(1) Let $g$ be a dg $L_{\infty}$-algebra and $\tau\in g^1$, we say that $\tau$ is a Maurer-Cartan element of $g$ if
\[
\sum_{k\geq 1} \frac{1}{k!} l_k(\tau,...,\tau)=0.
\]
The set of Maurer-Cartan elements of $g$ is noted $MC(g)$.

(2) The simplicial Maurer-Cartan set is then defined by
\[
MC_{\bullet}(g)=MC(g\hat{\otimes}\Omega_{\bullet}),
\],
where $\Omega_{\bullet}$ is the Sullivan cdga of de Rham polynomial forms on the standard simplex $\Delta^{\bullet}$ (see \cite{Sul})
and $\hat{\otimes}$ is the completed tensor product with respect to the filtration induced by $g$.
\end{defn}
The simplicial Maurer-Cartan set is a Kan complex, functorial in $g$ and preserves quasi-isomorphisms of complete $L_{\infty}$-algebras.
The Maurer-Cartan moduli set of $g$ is $\mathcal{MC}(g)=\pi_0MC_{\bullet}(g)$: it is the quotient of the set
of Maurer-Cartan elements of $g$ by the homotopy relation defined by the $1$-simplices.
When $g$ is a complete dg Lie algebra, it turns out that this homotopy relation is equivalent to the action of the gauge
group $exp(g^0)$ (a prounipotent algebraic group acting on Maurer-Cartan elements), so in this case
this moduli set coincides with the one usually known for Lie algebras.
We refer the reader to \cite{Yal4} for more details about all these results. The notion of Maurer-Cartan space allows to define the the classical deformation functor of $g$ given by
\begin{eqnarray*}
\underline{\mathcal{MC}}(g):Art_{\mathbb{K}} & \rightarrow & Set \\
R & \longmapsto & \mathcal{MC}(g\otimes m_R)
\end{eqnarray*}
and the derived deformation functor or derived formal moduli problem of $g$ given by
\begin{eqnarray*}
\underline{MC_{\bullet}}(g):dgArt_{\mathbb{K}}^{aug} & \rightarrow & sSet \\
R & \longmapsto & MC_{\bullet}(g\otimes m_R)
\end{eqnarray*}
(which belongs indeed to $FMP_{\mathbb{K}}$). By \cite[Corollary 2.4]{Yal4}, the tensor product $MC_{\bullet}(g\otimes m_R)$ does not need to be completed because $R$ is artinian. To see why Theorem 3.3 extends to $L_{\infty}$-algebras, let $\pi:L_{\infty}\stackrel{\sim}{\rightarrow}Lie$ be a cofibrant resolution of the operad $Lie$. This morphism induces a functor $p^*:dgLie\rightarrow L_{\infty}-Alg$ which associates to any dg Lie algebra the $L_{\infty}$-algebra with the same differential, the same bracket of arity $2$ and trivial higher brackets in arities greater than $2$. This functor fits in a Quillen equivalence
\[
p_{!}:L_{\infty}-Alg\leftrightarrows dgLie :p^*,
\]
where the left adjoint is a certain coequalizer (see \cite[Theorem 16.A]{Fre3}), and Quillen equivalences induce equivalences of the corresponding $\infty$-categories, so we have a commutative triangle of $\infty$-categories
\[
\xymatrix{
L_{\infty}-Alg\ar[dr]^-{\tilde{\psi}} & \\
dgLie \ar[u]^-{p^*}\ar[r]_-{\psi} & FMP_{\mathbb{K}}
}
\]
where $\psi$ and $\tilde{\psi}$ send a Lie algebra, respectively an $L_{\infty}$-algebra, to its derived formal moduli problem. The maps $p^*$ and $\psi$ are weak equivalences of $\infty$-categories, so $\tilde{\psi}:L_{\infty}-Alg\rightarrow FMP_{\mathbb{K}}$ is a weak equivalence of $\infty$-categories as well (here, by weak equivalence we mean a weak equivalence in the chosen model category of $\infty$-categories, say quasi-categories for instance).

\subsubsection{Twistings of $L_{\infty}$-algebras}

We recall briefly the notion of twisting by a Maurer-Cartan element. The twisting of a complete $L_{\infty}$ algebra $g$ by a Maurer-Cartan element $\tau$ is the complete $L_{\infty}$ algebra $g^{\tau}$
with the same underlying graded vector space and new brackets $l_k^{\tau}$ defined by
\[
l_k^{\tau}(x_1,...,x_k)=\sum_{i\geq 0}\frac{1}{i!}l_{k+i}(\underbrace{\tau,...,\tau}_i,x_1,...,x_k)
\]
where the $l_k$ are the brackets of $g$. The twisted $L_{\infty}$-algebra $g^{\varphi}$ is the \emph{deformation complex of $\varphi$}, that is, the derived formal moduli problem of $g^{\varphi}$ controls the deformation theory of $\varphi$. To see this, let us define another kind of Maurer-Cartan functor
\begin{eqnarray*}
\tilde{MC}_{\bullet}(g\otimes -):dgArt_{\mathbb{K}}^{aug} & \rightarrow & sSet \\
R & \longmapsto & MC_{\bullet}(g\otimes R).
\end{eqnarray*}
We replaced the maximal ideal $m_R$ in the definition of the deformation functor by the full algebra $R$. That is, the functor $\tilde{MC}_{\bullet}(g\otimes -)$ sends $R$ to the space of $R$-linear extensions of Maurer-Cartan elements of $g$. Then, for every augmented dg artinian algebra $R$ one has
\[
MC_{\bullet}(g^{\varphi}\otimes m_R)=hofib(MC_{\bullet}(g\otimes R)\rightarrow MC_{\bullet}(g),\varphi)
\]
where the map in the right side is induced by the augmentation $R\rightarrow\mathbb{K}$ and the homotopy fiber is taken over the base point $\varphi$. That is, the space $MC_{\bullet}(g^{\varphi}\otimes m_R)$ is the space of $R$-linear extensions of $\varphi$ as Maurer-Cartan elements of $g\otimes R$.

\section{Moduli spaces of algebraic structures}

\subsection{First version: a simplicial construction}

We refer the reader to \cite[Chapter 16, Chapter 17]{Hir} and \cite[]{Fre5} for some prerequisites about simplicial mapping spaces in model categories. We use this notion of simplicial mapping space and the model category structure on props to define our moduli spaces. Let us define a first version of this moduli space as a simplicial set. This was originally defined in the setting of simplicial operads \cite{Rez}, and can be extended to algebras over differential graded props as follows (see \cite{Yal6}):
\begin{defn}
Let $P_{\infty}$ be a cofibrant prop and $X$ be a cochain complex. The (simplicial) moduli space of $P_{\infty}$-algebra structures on $X$ is the simplicial set $P_{\infty}\{X\}$ defined
in each simplicial dimension $k$ by
\[
P_{\infty}\{X\}_k=Mor_{prop}(P_{\infty},End_X\otimes\Omega_k),
\]
where $(End_X\otimes\Omega_k)(m,n)=Hom(X^{\otimes m},X^{\otimes n})\otimes\Omega_k$.
\end{defn}
The Sullivan algebras $\Omega_k$ gather into a simplicial commutative differential graded algebra $\Omega_{\bullet}$ whose faces and degeneracies induce the simplicial structure on $P_{\infty}\{X\}$. The functor $(-)\otimes\Omega_{\bullet}$ is a functorial simplicial resolution in the model category of props \cite[Proposition 2.5]{Yal5}, so this simplicial moduli space is a homotopy mapping space in this model category. In particular, this means that this simplicial set is a Kan complex whose points are the $P_{\infty}$-algebra structures $P_{\infty}\rightarrow End_X$ and $1$-simplices are the homotopies between such structures (the prop $End_X\otimes\Omega_1$ forms a path object of $End_X$ in the model category of props). The later property implies that
\[
\pi_0P_{\infty}\{X\} = [P_{\infty},End_X]_{Ho(Prop)}
\]
is the set of homotopy classes of $P_{\infty}$-algebra structures on $X$. So our simplicial moduli space has the two first properties one expects from a moduli space: its points are the objects we want to classify and its connected components are the equivalence classes of these objects. Moreover, the fact that this is a homotopy mapping space implies that it is homotopy invariant with respect to the choice of a cofibrant resolution for the source, that is, any weak equivalence of cofibrant props $P_{\infty}\stackrel{\sim}{\rightarrow}Q_{\infty}$ induces a weak equivalence of Kan complexes
\[
Q_{\infty}\{X\}\stackrel{\sim}{\rightarrow}P_{\infty}\{X\}.
\]
So this is a well defined classifying object for homotopy $P$-algebra structures on $X$.

Another interesting homotopy invariant is the classification space of $P_{\infty}$-algebras, defined as the nerve $\mathcal{N}wP_{\infty}-Alg$ of the subcategory whose objects are $P_{\infty}$-algebras and morphisms are quasi-isomorphisms of $P_{\infty}$-algebras. By \cite{DK,DK2}, this classification space admits a decomposition
\[
\mathcal{N}wP_{\infty}-Alg \simeq \sqcap_{[X]\in\pi_0\mathcal{N}wP_{\infty}-Alg}\overline{W}L^HwP_{\infty}-Alg(X,X).
\]
Here the product ranges over weak equivalence classes of $P_{\infty}$-algebras, and $\overline{W}L^HwP_{\infty}-Alg(X,X)$ is the classifying complex of the simplicial monoid of zigzags of weak equivalences $X\stackrel{\sim}{\leftarrow}\bullet\stackrel{\sim}{\rightarrow}X$ in the hammock localization (or equivalently in the simplicial localization) of $P_{\infty}-Alg$ in the sense of Dwyer-Kan, i.e. the self equivalences of $X$ in the $\infty$-category of $P_{\infty}$-algebras. Let us note that when $P_{\infty}$ is an operad and $X$ is a cofibrant $P_{\infty}$-algebra, this space is equivalent to the usual simplicial monoid $haut_{P_{\infty}}(X)$ of self weak equivalences of $X$. This means that the classification space of $P_{\infty}$-algebras encodes symmetries and higher symmetries of $P_{\infty}$-algebras in their homotopy theory. Homotopy invariance of the classification space for algebras over props is a non trivial theorem:
\begin{thm}{\cite[Theorem 0.1]{Yal2}}
Let $\varphi:P_{\infty}\stackrel{\sim}{\rightarrow}Q_{\infty}$ be a weak equivalence between two cofibrant props. The map $\varphi$ gives rise to a functor
\[
\varphi^{*}: wQ_{\infty}-Alg\rightarrow wP_{\infty}-Alg\]

which induces a weak equivalence of simplicial sets
\[
\mathcal{N}\varphi^{*}:\mathcal{N}wQ_{\infty}-Alg\stackrel{\sim}{\rightarrow}\mathcal{N}wP_{\infty}-Alg.
\]
\end{thm}
Moreover, it turns out that the simplicial moduli space defined above gives a local approximation of this classification space, precisely we have the following result:
\begin{thm}{\cite[Theorem 0.1]{Yal6}}
Let $P_{\infty}$ be a cofibrant dg prop and $X$ be a cochain complex.
The commutative square
\[
\xymatrix{P_{\infty}\{X\}\ar[d]\ar[r] & \mathcal{N}wP_{\infty}-Alg\ar[d]\\
\{X\}\ar[r] & \mathcal{N}w\Ch
}
\]
is a homotopy pullback of simplicial sets.
\end{thm}
This homotopy fiber theorem has been applied to study the homotopy type of realization spaces in \cite{Yal7} in terms of derivation complexes and to count equivalence classes of realizations (of Poincaré duality for example).

The reader has probably noticed that we used the following property to define our simplicial moduli space: tensoring a prop by a cdga componentwise preserves the prop structure. This allows us to extend the definition of this moduli space and make it a simplicial presheaf of cdgas
\[
\underline{Map}(P_{\infty},Q):R\in CDGA_{\mathbb{K}}\mapsto Map_{Prop}(P_{\infty},Q\otimes A).
\]
Moreover, the notion of classification space defined above in the sense of Dwyer-Kan can also be extended to a simplicial presheaf. For this, we use that for any cdga $R$, the category $Mod_R$ is a (cofibrantly generated) symmetric monoidal model category tensored over chain complexes, so that one can define the category $P_{\infty}-Alg(Mod_R)$ of $P_{\infty}$-algebras in $Mod_R$. The assignment
\[
A\mapsto wP_{\infty}-Alg(Mod_R)
\]
defines a weak presheaf of categories in the sense of \cite[Definition I.56]{Ane}. It sends a morphism $A\rightarrow B$ to the symmetric monoidal functor $-\otimes_A B$ lifted at the level of $P_{\infty}$-algebras. This weak presheaf can be strictified into a presheaf of categories (see \cite[Section I.2.3.1]{Ane}). Applying the nerve functor then defines a simplicial presheaf of Dwyer-Kan classification spaces that we note $\underline{\mathcal{N}wP_{\infty}-Alg}$. The simplicial presheaf $\underline{\mathcal{N}wCh_{\mathbb{K}}}$ associated to $A\mapsto Mod_A$ is the simplicial presheaf of quasi-coherent modules of \cite[Definition 1.3.7.1]{TV}. The constructions above then make the following generalization of Theorem 4.3 meaningful:
\begin{prop}{\cite[Proposition 2.13]{GY}}
Let $P_{\infty}$ be a cofibrant prop and $X$ be a chain complex. The forgetful functor $P_{\infty}-Alg\rightarrow Ch_{\mathbb{K}}$ induces a homotopy fiber sequence
\[
\underline{P_{\infty}\{X\}}\rightarrow \underline{\mathcal{N}wP_{\infty}-Alg} \rightarrow \underline{\mathcal{N}wCh_{\mathbb{K}}}
\]
of simplicial presheaves over cdgas, taken over the base point $X$.
\end{prop}

\subsection{Second version: a stack construction and the associated deformation theory}

If $P$ is a properad with cofibrant resolution $(\mathcal{F}(s^{-1}C),\partial)\stackrel{\sim}{\rightarrow}P$ for a certain homotopy coproperad $C$ (see \cite[Section 4]{MV1} for the definition of homotopy coproperads), and $Q$ is any properad, then  we consider the total complex $g_{P,Q}=Hom_{\Sigma}(\overline{C},Q)$ given by homomorphisms of $\Sigma$-biobjects from the augmentation ideal of $C$ to $Q$. In the case $Q=End_X$ we will note it $g_{P,X}$. By \cite[Theorem 5]{MV2}, it is a complete dg $L_{\infty}$ algebra whose Maurer-Cartan elements are prop morphisms $P_{\infty}\rightarrow Q$. This $L_{\infty}$-structure was also independently found in \cite[Section 5]{Mar}, where it is proved that such a structure exists when replacing our cofibrant resolution above by the minimal model of a $\mathbb{K}$-linear prop (and its completeness follows by \cite[Proposition 15]{Mar}). In \cite{Yal5}, we prove a non trivial generalization of this result at the level of simplicial presheaves:
\begin{thm}{\cite[Theorem 2.10,Corollary 4.21]{Yal5}}
Let $P$ be a dg properad equipped with a minimal model $P_{\infty}:=(\mathcal{F}(s^{-1}C),\partial)\stackrel{\sim}{\rightarrow}P$ and $Q$ be a dg properad. Let us consider the simplicial presheaf
\[
\underline{Map}(P_{\infty},Q):R\in CDGA_{\mathbb{K}}\mapsto Map_{Prop}(P_{\infty},Q\otimes A)
\]
where $CDGA_{\mathbb{K}}$ is the category of commutative differential graded $\mathbb{K}$-algebras and $Q\otimes A$ is the componentwise tensor product defined by $(Q\otimes A)(m,n)=Q(m,n)\otimes A$.
This presheaf is equivalent to the simplicial presheaf
\[
\underline{\tilde{MC}_{\bullet}}(Hom_{\Sigma}(\overline{C},Q)):A\in CDGA_{\mathbb{K}}\mapsto MC_{\bullet}(Hom_{\Sigma}(\overline{C},Q)\otimes A)
\]
associated to the complete $L_{\infty}$-algebra $Hom_{\Sigma}(\overline{C},Q)$.
\end{thm}
In the case $Q=End_X$, we get the simplicial presheaf which associates to $A$ the moduli space of $P_{\infty}$-algebra structures on $X\otimes A$. Let us note that $\underline{Map}(P_{\infty},Q)$ can be alternately defined by
\[
A\mapsto Map_{Prop(Mod_A)}(P_{\infty}\otimes A,Q\otimes A),
\]
where $Map_{Prop(Mod_A)}$ is the simplicial mapping space in the category of props in dg $A$-modules.
In the case $Q=End_X$, we have $Q\otimes A\cong End_{X\otimes A}^{Mod_A}$ where $End_{X\otimes A}^{Mod_A}$ is the endormorphism prop of $X\otimes A$ taken in the category of $A$-modules. That is, it associates to $A$ the simplicial moduli space of $A$-linear $P_{\infty}$-algebra structures on $X\otimes A$ in the category of $A$-modules.
This theorem applies to a large class of algebraic structures, including for instance Frobenius algebras, Lie bialgebras and their variants such as involutive Lie bialgebras, as well as the properad $Bialg$ encoding associative and coassociative bialgebras.

Under additional assumptions, we can equip such a presheaf with a stack structure:
\begin{thm}{\cite[Corollary 0.8]{Yal5}}
(1) Let $P_{\infty}=(\mathcal{F}(s^{-1}C),\partial)\stackrel{\sim}{\rightarrow} P$ be a cofibrant resolution of a dg properad $P$ and $Q$ be any dg properad such that each $Q(m,n)$ is a bounded complex of finite dimension in each degree.
The functor
\[
\underline{Map}(P_{\infty},Q):A\in CDGA_{\mathbb{K}}\mapsto Map_{Prop}(P_{\infty},Q\otimes A)
\]
is an affine stack in the setting of complicial algebraic geometry of \cite{TV}.

(2) Let $P_{\infty}=(\mathcal{F}(s^{-1}C),\partial)\stackrel{\sim}{\rightarrow} P$ be a cofibrant resolution of a dg properad $P$ in non positively graded cochain complexes,
and $Q$ be any properad such that each $Q(m,n)$ is a finite dimensional vector space.
The functor
\[
\underline{Map}(P_{\infty},Q):A\in CDGA_{\mathbb{K}}\mapsto Map_{Prop}(P_{\infty},Q\otimes A)
\]
is an affine stack in the setting of derived algebraic geometry of \cite{TV}, that is, an affine derived scheme.
\end{thm}
In the derived algebraic geometry context, the derived stack $\underline{Map}(P_{\infty},Q)$ is not affine anymore wether the $Q(m,n)$ are not finite dimensional vector spaces. However, we expect these stacks to be derived $n$-Artin ind-stacks for the $Q(m,n)$ being perfect complexes with finite amplitude $n$, using the characterization of derived $n$-Artin stacks via resolutions by Artin $n$-hypergroupoids given in \cite{Pri}.

We denote by $\mathbb{T}_{\underline{Map}(P_{\infty},Q),x_{\varphi}}$ the tangent complex of $\underline{Map}(P_{\infty},Q)$ at an $A$-point $x_{\varphi}$ associated to a properad morphism $\varphi:P_{\infty}\rightarrow Q\otimes_e A$. As we explained before in Section 3, non-positive cohomology groups of the deformation complex correspond to negative groups of the tangent complex, which computes the higher automorphisms (higher symmetries) of the point, and the positive part which computes the obstruction theory. Adding some finiteness assumptions on the resolution $P_{\infty}$, we can explicit the ring of functions of this affine stack:
\begin{thm}{\cite[Theorem 0.14]{Yal5}}
Let $P$ be a dg properad equipped with a cofibrant resolution $P_{\infty}:=\Omega(C)\stackrel{\sim}{\rightarrow}P$,
where $C$ admits a presentation $C=\mathcal{F}(E)/(R)$, and $Q$ be a dg properad such that each $Q(m,n)$ is a bounded
complex of finite dimension in each degree. Let us suppose that each $E(m,n)$ is of finite dimension, and that there exists an integer $N$ such that $E(m,n)=0$ for $m+n>N$. Then

(1) The moduli stack $\underline{Map}(P_{\infty},Q)$ is isomorphic to $\mathbb{R}Spec_{C^*(Hom_{\Sigma}(\overline{C},Q))}$,
where $C^*(Hom_{\Sigma}(\overline{C},Q))$ is the Chevalley-Eilenberg algebra of $Hom_{\Sigma}(\overline{C},Q)$.

(2) The cohomology of the tangent dg Lie algebra at a $\mathbb{K}$-point $\varphi:P_{\infty}\rightarrow Q$
is explicitely determined by
\[
H^*(\mathbb{T}_{\underline{Map}(P_{\infty},Q),x_{\varphi}}[-1]) \cong
H^*(Hom_{\Sigma}(\overline{C},Q)^{\varphi}).
\]
\end{thm}
This theorem applies to a wide range of structures including for instance Frobenius algebras, Lie bialgebras and their variants such as involutive Lie bialgebras, and associative-coassociative bialgebras.

\subsection{Properties of the corresponding formal moduli problems and derived deformation theory}

Before turning to formal moduli problems, a natural question after reading the previous section is the following: how are the tangent complexes of our moduli spaces related to the usual cohomology theories of well-known sorts of algebras such as Hochschild cohomology of associative algebras, Harrison cohomology of commutative algebras, Chevalley-Eilenberg cohomology of Lie algebras, or Gerstenhaber-Schack cohomology of associative-coassociative bialgebras (introduced to study the deformation theory of quantum groups \cite{GS}). It turns out that these tangent Lie algebras do not give exactly the usual cohomology theories, but rather shifted truncations of them. For instance, let us consider the Hochschild complex $Hom (A^{\otimes >0}, A)$ of a dg associative algebra $A$.  This Hochschild complex is bigraded with a cohomological grading induced by the grading of $A$  and a weight grading given by the tensor powers $A^{\otimes \bullet}$. It turns out that the part $Hom(A,A)$ of weight $1$ in the Hochschild complex is the missing part in $g_{Ass,A}^{\varphi}$ (the $L_{\infty}$-algebra of Theorem 4.7, where $\varphi:Ass\rightarrow End_A$ is the associative algebra structure of $A$). There is also a ``full'' version of the Hochschild complex defined by $Hom (A^{\otimes \geq 0}, A)$. These three variants of Hochschild complexes give a sequence of inclusions of three dg Lie algebras
\[
Hom ( A^{\otimes \geq 0}, A)[1]\supset Hom ( A^{\otimes >0}, A)[1]\supset Hom ( A^{\otimes >1}, A)[1].
\]
All of these have been considered in various places in the litterature, but without comparison of their associated moduli problems. For the full complex, it is known that it controls the linear deformation theory of $Mod_A$ as a dg category \cite{KellerLowen, Preygel}.

The same kind of open question arises for other cohomology theories and their variants, and one of the achievements of our work with Gregory Ginot \cite{GY} was to describe precisely the moduli problems controlled by these variants and how they are related in the general context of $P_{\infty}$-algebras.

The formal moduli problem $\underline{P_{\infty}\{X\}}^{\varphi}$ controlling the formal deformations of a $P_{\infty}$-algebra structure $\varphi:P_{\infty}\rightarrow End_X$ on $X$ is defined, on any augmented dg artinian algebra $R$, by the homotopy fiber
\[
\underline{P_{\infty}\{X\}}^{\varphi}(R)=hofib(\underline{P_{\infty}\{X\}}(R)\rightarrow \underline{P_{\infty}\{X\}}(\mathbb{K}))
\]
taken over the base point $\varphi$, where the map is induced by the augmentation $R\rightarrow\mathbb{K}$.
The twisting of the complete $L_{\infty}$-algebra $Hom_{\Sigma}(\overline{C},End_X)$ by a properad morphism $\varphi:P_{\infty}\rightarrow End_X$ is the deformation complex of $\varphi$, and we have an isomorphism
\[
g_{P,X}^{\varphi} = Hom_{\Sigma}(\overline{C},End_X)^{\varphi} \cong Der_{\varphi}(\Omega(C),End_X)
\]
where the right-hand term is the complex of derivations with respect to $\varphi$ \cite[Theorem 12]{MV2}, whose $L_{\infty}$-structure induced by the twisting of the left-hand side is equivalent to the one of \cite[Theorem 1]{Mar}.
Section 3.2.2 combined with Theorem 4.5 tells us which formal moduli problem this deformation complex controls:
\begin{prop}{\cite[Proposition 2.11]{GY}}
The tangent $L_{\infty}$-algebra of the formal moduli problem $\underline{P_{\infty}\{X\}}^{\varphi}$ is given by
\[
g_{P,X}^{\varphi} = Hom_{\Sigma}(\overline{C},End_X)^{\varphi}.
\]
\end{prop}
In derived algebraic geometry, a Zariski open immersion of derived Artin stacks $F\hookrightarrow G$ induces a weak equivalence between the tangent complex over a given point of $F$ and the tangent complex over its image in $G$ \cite{TV}. It is thus natural to wonder more generally whether an ``immersion'' of an $\infty$-category $\mathcal{C}$ into another $\infty$-category $\mathcal{D}$ induces an equivalence between the deformation problem of an object $X$ of $\mathcal{C}$ (which should be in some sense a tangent space of $\mathcal{C}$ at $X$) and the deformation problem of its image in $\mathcal{D}$, in particular an equivalence of the corresponding tangent dg Lie algebras when such a notion makes sense. Here the word ``immersion'' has to be understood as ``fully faithful conservative $\infty$-functor'', that is, a fully faithful $\infty$-functor $\mathcal{C}\rightarrow\mathcal{D}$ such that a map of $\mathcal{C}$ is an equivalence if and only if its image in $\mathcal{D}$ is a weak equivalence. In the case of $\infty$-categories of algebras over props, Proposition 4.4 tells us that the formal moduli problem $\underline{P_{\infty}\{X\}}^{\varphi}$ is the ``tangent space'' over $(X,\varphi)$ to the Dwyer-Kan classification space of the $\infty$-category of $P_{\infty}$-algebras, with associated tangent $L_{\infty}$-algebra $g_{P,X}^{\varphi}$. In this setting, we can thus transform the intuition above into a precise statement:
\begin{thm}{\cite[Theorem 2.16]{GY}}
Let $F:P_{\infty}-Alg\rightarrow Q_{\infty}-Alg$ be a fully faithful and conservative $\infty$-functor inducing functorially in $A$, for every augmented artinian cdga $A$, a fully faithful and conservative $\infty$-functor $F:P_{\infty}-Alg(Mod_A)\rightarrow Q_{\infty}-Alg(Mod_A)$. Then $F$ induces an equivalence of formal moduli problems
\[
\underline{P_{\infty}\{X\}}^{\varphi} \sim  \underline{Q_{\infty}\{F(X)\}}^{F(\varphi)},
\]
where $F(\varphi)$ is the $Q_{\infty}$-algebra structure on the image $F(X,\varphi)$ of $X,\varphi$ under $F$, hence an equivalence of the associated $L_{\infty}$-algebras
\[
g_{P,X}^{\varphi}\sim g_{Q,F(X)}^{F(\varphi)}.
\]
\end{thm}

Proposition 4.4 also hints towards the fact that $g_{P,X}^{\varphi}$ does not control the deformation theory of homotopy automorphisms of $(X,\varphi)$ in the infinitesimal neighbourhood of $id_{(X,\varphi)}$, but should be closely related to it, since classification spaces decompose into disjoint unions of homotopy automorphisms. These homotopy automorphisms form a derived algebraic group \cite{Fra}, and as for underived algebraic groups, one can associate Lie algebras to such objects. Indeed, as explained in \cite{Fra}, given a moduli functor $F$ and a point $x\in F(\mathbb{K})$, the reduction of $F$ at $x$ is the functor $F_x$ defined by the homotopy fiber
\[
F_x(R)=hofib(F(R)\rightarrow F(\mathbb{K})
\]
where the map is induced by the augmentation of $R$ and the homotopy fiber is taken over the base point $x$. A point $x$ of $F(\mathbb{K})$ such that the reduction of $F$ at $x$ is a formal moduli problem (called an infinitesimal moduli problem in \cite[Definition 4.5]{Fra}) is called formally differentiable \cite[Definition 4.10]{Fra}, so there is a tangent Lie algebra of $F$ at $x$ defined as the Lie algebra of the formal moduli problem $F_x$. In the case of derived algebraic groups, the neutral element is a formally differentiable point and the Lie algebra of a derived algebraic group is the Lie algebra of its reduction at the neutral element. This is the natural extension to a derived framework of the well known Lie algebra of a Lie group. Consequently, there should be a homotopy fiber sequence of $L_{\infty}$-algebras relating $Lie(\underline{haut}_{P_{\infty}}(X,\varphi))$ to the tangent $L_{\infty}$-algebra $g_{P,X}^{\varphi}$ of $\underline{P_{\infty}\{X\}}^{\varphi}$.

Let us explicit a bit the construction of derived algebraic groups of homotopy automorphisms. Given a complex $X$, its homotopy automorphism group is noted $haut(X)$. Given a $P_{\infty}$-algebra $(X,\varphi)$, its homotopy automorphism group in the $\infty$-category of $P_{\infty}$-algebras is noted $haut_{P_{\infty}}(X,\varphi)$. In the general case, it is defined by Dwyer-Kan's hammock localization $L^HwP_{\infty}((X,\varphi),(X,\varphi))$, since we do not have a model category structure on the category of $P_{\infty}$-algebras. However in the particular case where $P_{\infty}$ is an operad, it turns out that this construction is equivalent to the usual simplicial monoid of homotopy automorphisms of $(X,\varphi)$ in the model category of $P_{\infty}$-algebras (the simplicial sub-monoid of self weak equivalences in the usual homotopy mapping space $Map_{P_{\infty}-Alg}(X,X)$).
\begin{rem}
What we mean here by a homotopy automorphism is a self weak equivalence, not the homotopy class of a strict automorphism.
\end{rem}
The derived algebraic group $\underline{haut}(X)$ of homotopy automorphisms of $X$ is defined by the strictification of the weak simplicial presheaf
\[
R\mapsto haut_{Mod_R}(X\otimes R),
\]
where $haut_{Mod_A}$ is the simplicial monoid of homotopy automorphisms in the category of $A$-modules. The derived algebraic group $\underline{haut}_{P_{\infty}}(X,\varphi)$ of homotopy automorphisms of $(X,\varphi)$ is defined by the strictification of the weak simplicial presheaf
\[
R\mapsto haut_{P_{\infty}}(X\otimes R,\varphi\otimes R)_{Mod_R}
\]
where $haut_{P_{\infty}}(X\otimes R,\varphi\otimes R)_{Mod_R}$ is the simplicial monoid of homotopy automorphisms
of $(X\otimes R,\varphi\otimes R)\in P_{\infty}-Alg(Mod_R)$. The reduction of $\underline{haut}_{P_{\infty}}(X,\varphi)$ at $id_{(X,\varphi)}$ associates to any augmented dg artinian algebra $R$ the space of $R$-linear extensions of homotopy automorphisms living in the connected component of $id_{(X,\varphi)}$, that is, homotopy isotopies. Finally, the deformation complex of $\varphi$ in the $\infty$-category of props and the deformation complex of homotopy isotopies of $(X,\varphi)$ in the $\infty$-category of $P_{\infty}$-algebras are related by the expected fiber sequence:
\begin{prop}{\cite[Proposition 2.14]{GY}}
There is a homotopy fiber sequence of $L_{\infty}$-algebras
\[
g_{P,X}^{\varphi}\rightarrow Lie(\underline{haut}_{P_{\infty}}(X,\varphi))\rightarrow Lie(\underline{haut}(X)).
\]
\end{prop}
Moreover, we can explicit $Lie(\underline{haut}_{P_{\infty}}(X,\varphi))$ as a slight modification $g_{P^+,X}^{\varphi^+}$ of $g_{P,X}^{\varphi}$, which consists in adding a component $Hom(X,X)$ to $g_{P,X}^{\varphi}$ (we refer to \cite[Section 3]{GY}):
\begin{thm}{\cite[Theorem 3.5]{GY}}
There is a quasi-isomorphism of $L_{\infty}$-algebras
\[
g_{P^+,X}^{\varphi^+}\simeq Lie(\underline{haut}_{P_{\infty}}(X,\varphi)).
\]
\end{thm} 
The conceptual explanation underlying this phenomenon is that $g_{P,X}^{\varphi}$ controls the deformations of the $P_{\infty}$-algebra structure over a fixed complex $X$, whereas $g_{P^+,X}^{\varphi^+}$ controls deformations of this  $P_{\infty}$-algebra structure plus compatible deformations of the differential of $X$, that is, deformations of the $P_{\infty}$-algebra structure up to self quasi-isomorphisms of $X$. This is the role of the part $Hom(X,X)$ appearing for instance in Hochschild cohomology. For instance, given an associative dg algebra $A$, the complex $g_{Ass^+,A}^{\varphi^+}\cong Hom ( A^{\otimes >0}, A) [1]$ 
computes the Hochschild cohomology of $A$ and the complex $g_{Ass,A}^{\varphi} \cong Hom ( A^{\otimes >1}, A)[1]$ is the one controlling the formal moduli problem of deformations of $A$ with fixed differential. The full shifted Hochschild complex $Hom ( A^{\otimes\geq 0}, A)[1)$ controls the linear deformations of the dg category $Mod_A$.

For an $n$-Poisson algebra $A$ (Poisson algebras with a Poisson bracket of degree $1-n$), we have the same kind of variants of $L_{\infty}$-algebras: the full shifted Poisson complex $CH_{Pois_n}(A)[n]$ \cite{CaWi}, the deformation complex $CH_{Pois_n}^{(\bullet>0)}(A)[n]$ introduced by Tamarkin \cite{Ta-deformationofd-algebra} which is the part of positive weight in the full Poisson complex, and the further truncation $CH_{Pois_n}^{(\bullet>1)}(A)[n]$. In \cite[Section 6]{GY}, we solve the open problem to determine which deformation problems these $L_{\infty}$-algebras control:
\begin{thm}
Let $A$ be an $n$-Poisson algebra.

(1) The truncation $CH_{Pois_n}^{(\bullet>1)}(A)[n]$ is the deformation complex $g_{Pois_n,R}^{\varphi}$ of the formal moduli problem $\underline{{Pois_{n}}_{\infty} \{A \}}^{\varphi}$ of homotopy $n$-Poisson algebra structures deforming $\varphi$.

(2) Tamarkin's deformation complex \emph{controls deformations of $A$ into dg-$Pois_n$-algebras}, that is, it is the tangent Lie algebra $g_{Pois_n^+,A}^{\varphi^+}$ of $\underline{haut_{Pois_n}}(A)$.
\end{thm}
\begin{rem}
We conjecture that the $L_{\infty}$-algebra structure of the full shifted Poisson complex $CH_{Pois_n}^*(A)[n]$ controls the deformations of $Mod_A$ into  $E_{n-1}$-monoidal dg categories. This should have interesting consequences for deformation quantization of $n$-shifted Poisson structures in derived algebraic geometry \cite{CPTVV,Toen-ICM}.
\end{rem}

\section{Gerstenhaber-Schack conjecture, Kontsevich formality conjecture and deformation quantization}

\subsection{$E_n$-operads, higher Hochschild cohomology and the Deligne conjecture}\label{S:EnHoch}

Recall that an $E_n$-operad is a dg operad quasi-isomorphic to the singular chains $C_*D_n$ of the little $n$-disks operad. We refer the reader to \cite[Volume I]{Fre5} for a comprehensive treatment of the construction and main properties of the little $n$-disks operads. These $E_n$-operads satisfy the following properties:
\begin{itemize}
\item[$\bullet$] There is an isomorphism $H_*E_1\cong Ass$, where $Ass$ is the operad of associative algebras.

\item[$\bullet$] For $n\geq 2$, there is an isomorphism $H_*E_n\cong Pois_n$ where $Pois_n$ is the operad of $n$-Poisson algebras.

\item[$\bullet$] For $n\geq 2$, the $E_n$-operads are formal, i.e. there is a quasi-isomorphism of operads
\[
E_n\stackrel{\sim}{\rightarrow}Pois_n.
\]
Modulo a technical assumption satisfied in particular by $C_*D_n$, this formality holds over $\mathbb{Q}$ \cite{FW}.
\end{itemize}
The formality of the little $n$-disks operad has a long story of intermediate formality results (for $n=2$ over $\mathbb{Q}$ in \cite{Tam2}, for $n\geq 2$ over $\mathbb{R}$ in \cite{Ko1,LaV}, finally an intrinsic formality result over $\mathbb{Q}$ in \cite{FW}). This formality is the key point to prove Deligne conjecture, which states the existence of a homotopy Gerstenhaber structure (that is, the $E_2$-algebra structure) of the Hochschild complex with product given by the usual cup product. This result provided in turn an alternative method for deformation quantization of Poisson manifolds \cite{Ko1,Ko2,Tam1,Tam2}.

The cohomology theory of $E_n$-algebras is called the higher Hochschild cohomology or $E_n$-Hochschild cohomology:
\begin{defn}
The (full) $E_n$-Hochschild complex of an $E_n$-algebra $A$ is the derived hom $CH^*_{E_n}(A,A)=\mathbb{R}Hom^{E_n}_A(A,A)$ in the category of (operadic) $A$-modules over $E_n$.
\end{defn}
Given an ordinary associative (or $E_1$) algebra $A$, the category of (operadic) $A$-modules over $E_1$ is the category of $A$-bimodules, so one recovers the usual Hochschild cohomology. Moreover, the aforementioned Deligne conjecture generalizes to $E_n$-algebras:
\begin{thm}(see \cite[Theorem 6.28]{GTZ} or \cite{Fra, Lur2})
The $E_n$-Hochschild complex $CH^*_{E_n}(A,A)$ of an $E_n$-algebra $A$ forms an $E_{n+1}$-algebra.
\end{thm}
The endomorphisms $Hom_{biMod_A}(A,A)$ of $A$ in the category $biMod_A$ of $A$-bimodules form nothing but the center $Z(A)$ of $A$. Deriving this hom object gives the Hochschild complex, and the Hochschild cohomology $A$ satisfies $HH^0(A,A)=Z(A)$. One says that the Hoschchild complex is the derived center of $A$, and the result above can then be reformulated  as \lq\lq{}the derived center of an $E_n$-algebra forms an $E_{n+1}$-algebra\rq\rq{}. This sentence has actually a precise meaning, because higher Hochschild cohomology can be alternately defined as a centralizer in the $\infty$-category of $E_n$-algebras. We refer the reader to \cite{Lur2} for more details about this construction.
Associated to an $E_n$-algebra $A$, one also has its cotangent complex $L_A$,  which classifies square-zero extensions of $A$ \cite{Fra, Lur2}, and its dual the tangent complex $T_A:= Hom^{E_n}_A(L_A,A)\cong \mathbb{R}Der(A,A)$.
\begin{thm}(see \cite{Fra,Lur2})
The shifted tangent complex $T_A[-n]$ of an $E_n$-algebra is an $E_{n+1}$-algebra, and is related to its $E_n$-Hochschild complex by a homotopy fiber sequence of $E_{n+1}$-algebras
\[
A[-n] \rightarrow T_A[-n]\rightarrow CH^*_{E_n}(A,A).
\]
\end{thm}

\subsection{From bialgebras to $E_2$-algebras}

A bialgebra is a complex equipped with an associative algebra structure and a coassociative coalgebra structure, such that the product is a coalgebra map and equivalently the coproduct is an algebra map. That is, bialgebras are equivalently algebras in coalgebras or coalgebras in algebras. Their cohomology theory is the Gerstenhaber-Schack cohomology \cite{GS}, intertwinning Hochschild cohomology of algebras and co-Hochschild cohomology of coalgebras.
Such structures naturally occur in algebraic topology (homology or cohomology of an $H$-space, for instance loop spaces), Lie theory (universal enveloping algebras, cohomology of Lie groups), representation theory (group rings, regular functions on algebraic groups, Tannaka-Krein duality), quantum field theory (renormalization Hopf algebras, AdS/CFT formalism)...Here we are going to focus on their prominent role in quantum group theory, \cite{Dri1, Dri2, EK1, EK2, GS, Mer1, Mer2}. As explained in Example 4, deformation quantization of Lie bialgebras produce quantum groups, whose categories of representations are particularly well behaved (modular tensor categories) and used to produce topological invariants via $3$-TQFTs \cite{RT}. It turns out that bialgebras are deeply related to $E_n$-algebras, via the natural occurrence of $E_n$-structures in deformation quantization and representation theory of quantum groups for instance, leading people to investigate the relationship between these two kinds of structures to understand various related problems on both sides. A hope in particular was to establish some equivalence between their respective deformation theories, maybe even their homotopy theories. One of the first goals of \cite{GY} was to embody this long-standing hope in a precise mathematical incarnation. A first crucial step is to relate bialgebras to a ``half-restricted'' kind of $E_2$-coalgebras by the following equivalence of $\infty$-categories: 
\begin{thm}{\cite[Theorem 0.1]{GY}}

(1) There exists a bar-cobar adjunction
\[
\mathcal{B}^{enh}_{E_1}: E_1-Alg^{0-con}(dgCog^{conil})\rightleftarrows E_1-Cog^{conil}(dgCog^{conil}):\Omega^{enh}_{E_1}
\]
inducing an equivalence of $\infty$-categories between nilpotent homotopy associative algebras in conilpotent dg coalgebras ($0$-connected conilpotent homotopy associative bialgebras) and conilpotent homotopy coassociative coalgebras in conilpotent dg coalgebras.

(2) The equivalence above induces an equivalence of $(\infty,1)$-categories
\[
E_1-Alg^{aug,nil}(dgCog^{conil})\rightleftarrows E_1-Cog^{conil,pt}(dgCog^{conil})
\]
between nilpotent augmented conilpotent homotopy associative bialgebras and pointed conilpotent homotopy coassociative coalgebras in conilpotent dg coalgebras.
\end{thm}
In part (1), the notation $0-con$ means $0$-connected, that is dg bialgebras concentrated in positive degrees. In part (2), the notation $aug,nil$ stands for augmented and nilpotent. A typical example of such a bialgebra is the total complex of the symmetric algebra over a cochain complex. The notation $pt$ stands for pointed coalgebras, that is, a coalgebra $C$ equipped with a counit $\epsilon:C\rightarrow\mathbb{K}$ and a coaugmentation $e:\mathbb{K}\rightarrow C$ such that $\epsilon\circ e = id_{\mathbb{K}}$. More generally, one can wonder, working in a given stable symmetric monoidal $\infty$-category (not necessarily cochain complexes), under which conditions a bar-cobar adjunction induces an equivalence of $\infty$-categories between algebras over an operad and conilpotent coalgebras over its bar construction \cite{FG}.
Theorem 6.3 solves this conjecture of Francis-Gaitsgory \cite{FG} in the case where the base category is the category of conilpotent dg coalgebras, respectively the category of pointed conilpotent dg coalgebras.
\begin{rem}
This example is also interesting with respect to the conditions imposed on the ground symmetric monoidal $\infty$-category in \cite{FG}, since the categories considered here are a priori not pronilpotent in the sense of \cite[Definition 4.1.1]{FG}.
\end{rem}
Using Koszul duality of $E_n$-operads and an $\infty$-categorical version of Dunn's theorem \cite{Lur2,GTZ1}, we deduce from these equivalences the precise and long awaited relationship between homotopy theories of bialgebras and $E_2$-algebras.
The correct answer to this problem needs an appropriate notion of ``cobar construction for bialgebras'', which intertwines a bar construction on the algebra part of the structure with a cobar construction on the resulting $E_2$-coalgebra:
\begin{thm}{\cite[Corollary 0.2]{GY}}
The left adjoint of Theorem 0.1(2) induces a conservative fully faithful $\infty$-functor
\[
\tilde{\Omega}: E_1-Alg^{aug,nil}(dgCog^{conil})\hookrightarrow E_2-Alg^{aug}
\]
embedding augmented nilpotent and conilpotent homotopy associative bialgebras into augmented $E_2$-algebras.
\end{thm}
By Theorem 4.12, this ``immersion'' of $\infty$-categories induces equivalences of formal moduli problems between the moduli problem of homotopy bialgebra structures on a bialgebra $B$ and the moduli problem of $E_2$-algebra structures on its cobar construction $\tilde{\Omega} B$. Moreover, at the level of formal moduli problems controlling homotopy isotopies of $B$ and homotopy isotopies of $\tilde{\Omega}B$, the tangent $L_{\infty}$-algebras can be identified respectively with the shifted Gerstenhaber-Schack complex of $B$ and the shifted (truncated) higher Hochschild complex (or $E_2$-tangent complex) of $\Omega B$ as $L_{\infty}$-algebras:
\begin{thm}{\cite[Theorem 0.6]{GY}}
Let $B$ be a pointed conilpotent homotopy associative dg bialgebra. Let $\varphi:Bialg_{\infty}\rightarrow End_B$ be this homotopy bialgebra structure on $B$ (where $Bialg$ is the prop of associative-coassociative bialgebras) , and let $\tilde{\Omega}\varphi:E_2\rightarrow End_{\tilde{\Omega}B}$ be the corresponding $E_2$-algebra structure on its cobar construction $\tilde{\Omega} B$.

(1) There is a homotopy equivalence of formal moduli problems
\[
\underline{Bialg_{\infty}\{B\}}^{\varphi} \simeq \underline{E_2\{\tilde{\Omega}B\}}^{\tilde{\Omega}\varphi}.
\]
This homotopy equivalence induces a quasi-isomorphism of $L_{\infty}$-algebras
\[
g_{Bialg,B}^{\varphi}\stackrel{\sim}{\rightarrow} g_{E_2,\tilde{\Omega}B}^{\tilde{\Omega}\varphi}.
\]

(2) There is a homotopy equivalence of formal moduli problems
\[
\underline{Bialg_{\infty}^+}\{B\}^{\varphi^+} \simeq \underline{E_2^+}\{\tilde{\Omega}B\}^{\tilde{\Omega}\varphi^+}.
\]
This homotopy equivalence induces a quasi-isomorphism of $L_{\infty}$-algebras
\[
C^*_{GS}(B,B)[2]\stackrel{\sim}{\rightarrow} T_{\tilde{\Omega}(B)} 
\]
between the shifted Gerstenhaber-Schack complex of $B$ and the (truncated) $E_2$-Hochschild complex or $E_2$-tangent complex of $\tilde{\Omega}(B)$.
\end{thm}
The $L_{\infty}$ structure on the $E_2$-Hochschild complex $T_{\tilde{\Omega}(B)}$ is the one induced by the $E_3$ structure on $T_{\tilde{\Omega}(B)}[-2]$ (see Theorem 5.3). Proving that the higher Hochschild complex of the cobar construction of a bialgebra is a deformation complex of this bialgebra is important, since it allows to reduce questions of deformations of bialgebras to those of $E_2$-structures for which more tools are available.

\subsection{Gerstenhaber-Schack conjecture}

At the beginning of the 90's, Gerstenhaber and Schack enunciated (in a wrong way) a conjecture \cite{GS} characterizing the structure of the complex controlling the deformation theory of bialgebras,  which remained quite mysterious for a while. It is a dg bialgebra version of the Deligne conjecture. In \cite[Section 8]{GS}, we extended the equivalences of Theorem 5.7 to an equivalence of homotopy fiber sequences of $E_3$-algebras, getting a much stronger version of the longstanding Gerstenhaber-Schack conjecture for the different versions of the Gerstenhaber-Schack and $E_2$-Hochschild complexes:
\begin{thm}(Generalized Gerstenhaber-Schack conjecture \cite[Corollary 0.7]{GY})

(1) There is an $E_3$-algebra structures  on $C^*_{GS}(B,B)$ and a unital $E_3$-algebra structure on $C_{GS}^{full}(B,B))$
such that the following diagram  \[\xymatrix{ \tilde{\Omega}B [-1] \ar[r] & T_{\tilde{\Omega}(B)} \ar[r] &  CH_{E_2}^{*}( \tilde{\Omega}B, \tilde{\Omega}B) \\ 
  \tilde{\Omega}B[-1] \ar@{=}[u] \ar[r] & C_{GS}^{*}(B,B) \ar[u]^{\simeq} \ar[r] & C_{GS}^{full}(B,B) \ar[u]_{\simeq} }\] 
  is 
a commutative diagram of non-unital $E_3$-algebras with vertical arrows being equivalences. 

(2) The $E_3$-algebra structure on $C^*_{GS}(B,B)$ is a refinement of its $L_{\infty}$-algebra structure controlling the deformation theory of the bialgebra $B$.
\end{thm}
Let us note that the upper fiber sequence of part (1) is the fiber sequence of Theorem 5.3. In particular, the $E_3$-algebra structure on the deformation complex of dg bialgebra $B$ comes from the $E_3$-algebra structure on the $E_2$-Hochschild complex of $\tilde{\Omega}B$ given by the higher Deligne conjecture.
 
\subsection{Kontsevich formality conjecture and deformation quantization of Lie bialgebras}

Let us first recall briefly how deformation quantization of Poisson manifolds works in Kontsevich's work \cite{Ko2}. We fix a a finite dimensional Poisson manifold $M$, and we consider two complexes one can associate to such a manifold. First, the Hochschild complex $CH^*(\mathcal{C}^{\infty}(M),\mathcal{C}^{\infty}(M))$, second, the complex of polyvector fields $T_{poly}(M)=\left( \bigoplus_{k\geq 0}\bigwedge^k\Gamma T(M)[-k] \right)[1]$ where $\Gamma T(M)$ is the space of sections of the tangent bundle on $M$. The Poisson structure we fixed on $M$ is the datum of a bivector $\Pi\in\bigwedge^2\Gamma T(M)$ satifying the Maurer-Cartan equation, that is, a Maurer-Cartan element of weight $2$ in the Lie algebra of polyvector fields $T_{poly}(M)[1]$ (equipped with the Schouten-Nihenjuis bracket). To get the equivalent definition of Poisson manifold as a manifold whose ring of functions is a Poisson algebra, set $\{f,g\}=\Pi(df,dg)$. A well known theorem called the Hochschild-Kostant-Rosenberg theorem (HKR for short) states that the cohomology of $CH^*(\mathcal{C}^{\infty}(M),\mathcal{C}^{\infty}(M))$ is precisely $T_{poly}(M)$. In \cite{Ko2}, Kontsevich proved that there exists a $L_{\infty}$-quasi-isomorphism
\[
T_{poly}(M)[1] \stackrel{\sim}{\rightarrow} CH^*(\mathcal{C}^{\infty}(M),\mathcal{C}^{\infty}(M))[1]
\]
realizing in particular the isomorphism of the HKR theorem. We did not use the notion of $L_{\infty}$-quasi-isomorphism before, let us just say briefly that it is a quasi-isomorphism of cdgas between the Chevalley-Eilenberg algebra of $T_{poly}(M)[1]$ and the Chevalley-Eilenberg algebra of $CH^*(\mathcal{C}^{\infty}(M),\mathcal{C}^{\infty}(M))[1]$. In particular, it is determined by an infinite collection of maps $T_{poly}(M)[1]\rightarrow \Lambda^k(CH^*(\mathcal{C}^{\infty}(M),\mathcal{C}^{\infty}(M))[1])$ for $k\in\mathbb{N}$, whose first map is the HKR quasi-isomorphism.
\begin{rem}
An $L_{\infty}$-quasi-isomorphism of dg Lie algebras is actually equivalent to a chain of
quasi-isomorphisms of dg Lie algebras.
\end{rem}
This formality theorem then implies the deformation quantization of Poisson manifolds by the following arguments. First, noting $g[[\hbar]]_+=\bigoplus\hbar g^n[[\hbar]]$, one proves that the Maurer-Cartan set $MC(T_{poly}(M)[[\hbar]]_+)$ is the set of Poisson algebra structures on $\mathcal{C}^{\infty}(M)[[\hbar]]$ and that the Maurer-Cartan set $MC(D_{poly}(M)[[\hbar]]_+)$ is the set of $*_{\hbar}$-products, which are assocative products on $\mathcal{C}^{\infty}(M)[[\hbar]]$ of the form $a.b + B_1(a,b)t+...$ (i.e these products restrict to the usual commutative associative product on $\mathcal{C}^{\infty}(M)$).
Second, an $L_{\infty}$-quasi-isomorphism of nilpotent dg Lie algebras induces a bijection between the corresponding moduli sets of Maurer-Cartan elements, so there is a one-to-one correspondence between gauge equivalence classes of both sides. Consequently, isomorphism classes of formal Poisson structures on $M$ are in bijection with equivalence classes of $*_{\hbar}$-products.

Kontsevich builds explicit formality morphisms in the affine case $M=\mathbb{R}^d$, with formulae involving integrals
on compactification of configuration spaces and deeply related to the theory of multi-zeta functions. An alternative proof of the formality theorem for $M=\mathbb{R}^d$ due to Tamarkin \cite{Tam1}, relies on the formality of $E_2$-operads (hence on the choice of a Drinfeld associator) and provides a formality quasi-isomorphism of homotopy Gerstenhaber algebras (that is $E_2$-algebras)
\[
T_{poly}(\mathbb{R}^n) \stackrel{\sim}{\rightarrow} CH^*(\mathcal{C}^{\infty}(\mathbb{R}^n),\mathcal{C}^{\infty}(\mathbb{R}^n)).
\]
His method works as follows:
\begin{itemize}
\item[$\bullet$] Prove the Deligne conjecture stating the existence of an $E_2$-algebra structure on the Hochschild complex;

\item[$\bullet$] Transfer this structure along the HKR quasi-isomorphism to get an $E_2$-quasi-isomorphism betweeen $CH^*(\mathcal{C}^{\infty}(\mathbb{R}^n),\mathcal{C}^{\infty}(\mathbb{R}^n))$ with its $E_2$-algebra structure coming from the Deligne conjecture, and $T_{poly}(\mathbb{R}^n)$.

\item[$\bullet$] By the formality of $E_2$, this means that we have two $E_2$-structures on $T_{poly}(\mathbb{R}^n)$, the one coming from the Deligne conjecture and the one coming from the $Pois_2$-structure given by the wedge product and the Schouten-Nijenhuis bracket. One proves that $T_{poly}(M)$ has a unique homotopy class of $E_2$-algebra structures by checking that it is intrinsically formal (precisely, the $Aff(\mathbb{R}^n)$-equivariant $Pois_2$-cohomology of $T_{poly}(\mathbb{R}^n)$ is trivial).
\end{itemize}
This ``local'' formality for affine spaces is then globalized to the case of a general Poisson manifold by means of formal geometry \cite[Section 7]{Ko2}.

In the introduction of his celebrated work on deformation quantization of Poisson manifolds \cite{Ko2}, Kontsevich conjectured that a similar picture should underline the deformation quantization of Lie bialgebras. Etingof-Kazhdan quantization (see \cite{Dri2, EK1, EK2}) should be the consequence of a deeper formality theorem for the deformation complex $Def(Sym(V))$ of the symmetric bialgebra $Sym(V)$ on a vector space $V$. This deformation complex should possess an $E_3$-algebra structure whose underlying $L_{\infty}$-structure controls the deformations  of $Sym(V)$, and should be formal as an $E_3$-algebra. Then this formality result should imply a one-to-one correspondence between gauge classes of Lie bialgebra structures on $V$ and gauge classes of their quantizations. In \cite{GY}, we solved this longstanding conjecture at a greater level of generality than the original statement, and deduced a generalization of Etingof-Kadhan's deformation quantization theorem. Here we consider not a vector space, but a $\mathbb{Z}$-graded cochain complex $V$ whose cohomology is of finite dimension in each degree. By the results explained in Section 5.3, we know the existence of such an $E_3$-algebra structure (interestingly coming from the higher Hochschild complex of $\tilde{\Omega}Sym(V)$). It remains to prove the $E_3$-formality of $Def(Sym(V))$ by proving the homotopy equivalence between two $E_3$-structures on the Gerstenhaber-Schack cohomology of $Sym(V)$: the one transferred from $Def(Sym(V))$, and the canonical one coming from the action of $Pois_3$ (giving an $E_3$-structure via the formality $E_3\stackrel{\sim}{\rightarrow}Pois_3$). Indeed, the cohomology of $Def(Sym(V))$ (which is precisely the Gerstenhaber-Schack complex of $V$) is explicitely computable, and given by
\[
H^*_{GS}(Sym(V),Sym(V))\cong \hat{Sym}(H^*V[-1]\oplus H^*V^{\vee}[-1])
\]
where $\hat{Sym}$ is the completed symmetric algebra and $H^*V^{\vee}$ is the dual of $H^*V$ as a graded vector space. This symmetric algebra has a canonical $3$-Poisson algebra structure induced by the evaluation pairing between $H^*V$ and $H^*V^{\vee}$. In the spirit of Tamarkin's method, we have to use obstruction theoretic methods to show that $Def(Sym(V))$ is rigid as an $E_3$-algebra. We thus get a generalization of Kontsevich's conjecture (originally formulated in the case where $V$ is a vector space):
\begin{thm}(Kontsevich formality conjecture \cite[Theorem 0.8]{GY})
The deformation complex of the symmetric bialgebra $Sym(V)$ on a $\mathbb{Z}$-graded cochain complex $V$ whose cohomology is of finite dimension in each degree is formal over $\mathbb{Q}$ as an $E_3$-algebra.
\end{thm}
We prove it by using in particular the relationship between Gerstenhaber-Schack cohomology and $E_2$-Hochschild cohomology and the higher HKR-theorem for the latter \cite{CaWi}. We then obtain a new proof of Etingof-Kazhdan quantization theorem from the underlying $L_{\infty}$-formality given by our $E_3$-formality. Indeed, this formality induces an equivalence of the associated derived formal moduli problems, in particular we have an equivalence of Maurer-Cartan moduli sets (suitably extended over formal power series in one variable). On the right hand side, the Maurer-Cartan moduli set is identified with equivalence classes of homotopy Lie bialgebra structures on the cochain complex $V[[\hbar]]$. On the right hand side, it is identified with deformation quantization of these Lie bialgebras (formal deformations of $Sym(V)$ as a homotopy dg bialgebra). Moreover, what we get is actually a generalization of Etingof-Kazhdan quantization to homotopy dg Lie bialgebras:
\begin{cor}{\cite[Corollary 0.9]{GY}}
The $L_{\infty}$-formality underlying Theorem 5.9 induces a generalization of Etingof-Kazdhan deformation quantization theorem to homotopy dg Lie bialgebras whose cohomology is of finite dimension in each degree. In the case where $V$ is a vector space, this gives a new proof of Etingof-Kazdhan's theorem.
\end{cor}
This result encompasses the case of usual Lie bialgebras, because if $V$ is concentrated in degree $0$, then homotopy Lie bialgebra structures on $V$ are exactly Lie bialgebra structures on $V$.
\begin{rem}
Actually, what we prove in \cite{GY} is even stronger. We get a sequence of $E_3$-formality morphisms for the three variants of the Gerstenhaber-Shack complex \cite[Theorem 7.2]{GY}, indicating that important variants of deformation quantization like \cite{EH1} should also follow from such $E_3$-formality morphisms.
\end{rem}


\begin{thebibliography}{10}

\bibitem{Ane}M. Anel, \textit{Champs de modules des catégories linéaires et abéliennes}, PhD thesis, Université
Toulouse III – Paul Sabatier (2006).

\bibitem{AF}D. Ayala, J. Francis, \textit{Poincaré/Koszul duality}, preprint arXiv:1409.2478.

\bibitem{Baues} H.-J. Baues, \textit{The cobar construction as a Hopf algebra}, Invent. Math. {\bf 132} (1998), no.~3, 467--489.

\bibitem{BF}C. Berger, B. Fresse, \textit{Combinatorial operad actions on cochains}, Math. Proc. Cambridge Philos. Soc. 137 (2004), 135–174.

\bibitem{Ber}A. Berglund, \textit{Rational homotopy theory of mapping spaces via Lie theory for $L_{\infty}$-algebras}, 
Homology, Homotopy and Applications 17 (2015), 343-369.


\bibitem{BDG}D. Blanc, W. G. Dwyer, P. G. Goerss, \textit{The realization space of a $\Pi$-algebra: a moduli problem in algebraic topology}, Topology 43 (2004), 857-892.

\bibitem{CaWi} D. Calaque, T. Willwacher, \textit{Triviality of the higher Formality Theorem},  Proc. Amer. Math. Soc. 143 (2015), no. 12, 5181–5193.

\bibitem{CPTVV}D. Calaque, T. Pantev, B. Toen, M. Vaquié, G. Vezzosi, \textit{Shifted Poisson structures and deformation quantization}, preprint arXiv:1506.03699.

\bibitem{CMW}R. Campos, S. Merkulov, T. Willwacher, \textit{The Frobenius properad is Koszul}, Duke Math. Journal 165 (2016), no. 15, 2921-2989.

\bibitem{CS1}M. Chas, D. Sullivan, \textit{String Topology}, preprint arXiv:math/9911159.

\bibitem{CS2}M. Chas, D. Sullivan, \textit{Closed string operators in topology leading to Lie bialgebras and higher string algebra}, The legacy of Niels Henrik Abel, 771–784, Springer, Berlin, 2004.

\bibitem{CL} K. Cieliebak, J. Latschev, \textit{The role of string topology in symplectic field theory}, New perspectives and challenges in symplectic field theory, CRM Proc. Lecture Notes 49, Amer. Math. Soc. (2009), 113-146.

\bibitem{CFL}K. Cieliebak, K. Fukaya, J. Latschev, \textit{Homological algebra related to surfaces with boundary}, preprint arXiv:1508.02741 (2015).

\bibitem{Dri1} V. G. Drinfeld, Quantum groups, in {\it Proceedings of the International Congress of Mathematicians, 
Vol. 1, 2 (Berkeley, Calif., 1986)}, 798--820, Amer. Math. Soc., Providence, RI.

\bibitem{Dri2}V. G. Drinfeld, \textit{On quasitriangular quasi-Hopf algebras and on a group that is closely connected with $Gal(\overline{Q}/Q$}, Algebra i Analiz 2 (1990), 149-181. English translation in Leningrad Math. J. 2 (1991), 829-860.

\bibitem{DK}W. G. Dwyer, D. Kan, \textit{Function complexes in homotopical algebra}, Topology 19 (1980), 427-440.

\bibitem{DK2}W. G. Dwyer, D. Kan, \textit{A classification theorem for diagrams of simplicial sets}, Topology 23 (1984), 139-155.

\bibitem{EE}B. Enriquez, P. Etingof, \textit{On the invertibility of quantization functors}, J. Algebra 289 (2005), no. 2, 321-345.

\bibitem{EH1}B. Enriquez, G. Halbout, \textit{Quantization of quasi-Lie bialgebras},  J. Amer. Math. Soc. 23 (2010), 611–653. 

\bibitem{EH2}B. Enriquez, G. Halbout, \textit{Quantization of coboundary Lie bialgebras},  Ann. of Math. (2) 171 (2010), 1267–1345. 

\bibitem{EK1}P. Etingof, D. Kazhdan, \textit{Quantization of Lie bialgebras I}, Selecta Math. (N. S.) 2 (1996), no. 1, 1-41.

\bibitem{EK2}P. Etingof, D. Kazdhan, \textit{Quantization of Lie bialgebras. II, III}, Selecta Math. (N.S.) 4
    (1998), no. 2, 213–231, 233–269.
    
\bibitem{FT}Y. Félix, J.-C. Thomas, \textit{Rational BV-algebra in string topology}, Bull. Soc. Math. France 136 (2008), 311-327.

\bibitem{FFRS}J. Fjelstad, J. Fuchs,I. Runkel, C. Schweigert, \textit{Topological and conformal field theory as Frobenius algebras}, Categories in algebra, geometry and mathematical physics, 225–247, Contemp. Math., 431, Amer. Math. Soc., Providence, RI, 2007.

\bibitem{FG}J. Francis, D. Gaitsgory, \textit{Chiral Koszul duality}, Selecta Math. New Ser. 18 (2012), 27-87.

\bibitem{Fra}J. Francis, \textit{The tangent complex and Hochschild cohomology of $E_n$-rings}, Compositio Mathematica 149 (2013), no. 3, 430-480.

\bibitem{Free}D. S. Freed, \textit{Remarks on Chern-Simons theory}, Bulletin of the AMS 46 (2009), 221-254.

\bibitem{Fre3}B. Fresse, \textit{Modules over operads and functors}, Lecture Notes in Mathematics 1967, Springer-Verlag (2009).

\bibitem{Frep}B. Fresse, \textit{Props in model categories and homotopy invariance of structures}, Georgian Math. J. 17 (2010), pp. 79-160.

\bibitem{Fre4}B. Fresse, \textit{Koszul duality of En-operads}, Selecta Math. (N.S.) 17 (2011), 363–434.

\bibitem{Fre5}B. Fresse, \textit{Homotopy of operads and Grothendieck-Teichmuller groups I and II}, to appear in  Mathematical Surveys and Monographs, AMS (available at http://math.univ-lille1.fr/~fresse/OperadHomotopyBook/).

\bibitem{FW}B. Fresse, T. Willwacher, \textit{The intrinsic formality of $E_n$-operads}, preprint arXiv:1503.08699.

\bibitem{GTV}I. Galvez-Carillo, A. Tonks, B. Vallette, \textit{Homotopy Batalin-Vilkovisky algebras}, J. Noncommut. Geom. 6 (2012), 539-602.

\bibitem{GS}M. Gerstenhaber, S. D. Schack, \textit{Algebras, bialgebras, quantum groups, and algebraic deformations}, Deformation theory and quantum groups with applications to mathematical physics (Amherst, MA, 1990), Contemp. Math. 134 (1992), 51-92.

\bibitem{GJ}E. Getzler, J. D. S. Jones \textit{Operads, homotopy algebra and iterated integrals for double loop spaces}, preprint arXiv:hep-th/9403055 (1994).

\bibitem{GiS}J. Giansiracusa, P. Salvatore, \textit{Formality of the framed little 2-discs operad and semidirect products}, Homotopy theory of function spaces and related topics, Contemp. Math. 519 (2010), 115–121.
   
\bibitem{GH}G. Ginot, G. Halbout, \textit{A Formality Theorem for Poisson Manifolds}, Lett. Math. Phys. 66 (2003), 37-64.

\bibitem{GTZ1}  G. Ginot, T. Tradler, M. Zeinalian, \textit{Higher Hochschild Homology, Topological Chiral
Homology and Factorization Algebras}, Comm. Math. Phys. 326 (2014), no. 3, 635-686.

\bibitem{GTZ}G. Ginot, T. Tradler, M. Zeinalian, \textit{Higher Hochschild cohomology of E-infinity algebras, Brane topology and centralizers of E-n algebra maps}, preprint arXiv:1205.7056.

\bibitem{GoH}P. G. Goerss, M. J. Hopkins, \textit{Moduli spaces of commutative ring spectra}, Structured ring spectra, London Math. Soc. Lecture Note Ser. 315 (2004), 151-200.

\bibitem{Gol}W. M. Goldman, \textit{The symplectic nature of fundamental groups of surfaces}, Adv. in Math. 54 (1984), 200-225.

\bibitem{Gol2}W. M. Goldman, \textit{Invariant functions on Lie groups and Hamiltonian flows of surface group representations}, Invent. Math. 85 (1986), no. 2, 263–302.

\bibitem{GM}W. M. Goldman, J. Millson, \textit{The deformation theory of representations of fundamental groups of compact Kähler manifolds}, Publ. Math. IHES 67 (1988), 43-96.

\bibitem{Hen}B. Hennion, \textit{Tangent Lie algebra of derived Artin stacks}, to appear in J. für die reine und angewandte Math. (Crelles Journal).

\bibitem{Hin0}V. Hinich, \textit{DG coalgebras as formal stacks}, J. Pure Appl. Algebra 162 (2001), 209-250.

\bibitem{Hin}V. Hinich, \textit{Tamarkin’s proof of Kontsevich formality theorem}, Forum Math. 15 (2003), 591–614.

\bibitem{Hov}Mark Hovey, \textit{Model categories}, Mathematical Surveys and Monographs volume 63, AMS (1999).

\bibitem{Hir}Philip S. Hirschhorn, \textit{Model categories and their localizations}, Mathematical Surveys and Monographs volume 99, AMS (2003).

\bibitem{Kad}T. Kadeishvili, \textit{On the cobar construction of a bialgebra}, Homology Homotopy Appl. 7 (2005), 109–122.

\bibitem{Kap-TFT} A. Kapustin, \textit{Topological field theory, higher categories, and their applications},
in {\it Proceedings of the International Congress of Mathematicians. Volume III}, Hindustan Book Agency, New Delhi (2010), 2021--2043.

\bibitem{Kel}B. Keller, \textit{A-infinity algebras, modules and functor categories}, in Trends in representation theory of  algebras and related topics,  Contemp. Math. 406 (2006), 67-93.

\bibitem{KellerLowen} B. Keller\ and\ W. Lowen, \textit{On Hochschild cohomology and Morita deformations},
Int. Math. Res. Not. IMRN {\bf 2009}, no.~17, 3221--3235.

\bibitem{Koc}J. Kock, \textit{Frobenius algebras and 2D topological quantum field theories},
London Mathematical Society Student Texts No. 59, 2003.

\bibitem{Ko1}M. Kontsevich, \textit{Operads and motives in deformation quantization}, Moshé Flato (1937-1998), Lett. Math. Phys. 48 (1999), no. 1, 35-72.

\bibitem{Ko2}M. Kontsevich, \textit{Deformation quantization of Poisson manifolds}, Lett. Math. Phys. 66 (2003), no. 3, 157-216.

\bibitem{LaV}P. Lambrechts,  and I. Volic, \textit{Formality of the little N-disks operad}, Mem. Amer. Math. Soc. 230 (1079), 2014.

\bibitem{LP}A. D. Lauda, H. Pfeiffer, \textit{Open-closed strings: two-dimensional extended TQFTs and Frobenius algebras}, Topology and its applications 155 (2008), 623-666.

\bibitem{LV}J-L. Loday, B. Vallette, \textit{Algebraic Operads}, Grundlehren der mathematischen Wissenschaften,
Volume 346, Springer-Verlag (2012).

\bibitem{Lur0}J. Lurie, \textit{Derived Algebraic Geometry X}, available at http://www.math.harvard.edu/~lurie/

\bibitem{Lur1}J. Lurie, \textit{Higher topos theory}, Annals of Mathematics Studies 170, Princeton University Press, Princeton, NJ, 2009.

\bibitem{Lur2}J. Lurie, \textit{Higher Algebra}, book available at http://www.math.harvard.edu/~lurie/

\bibitem{MLa}S. MacLane, \textit{Categorical algebra}, Bull. Amer. Math. Soc. 71 (1965), 40-106.

\bibitem{Man}M. A. Mandell, \textit{Cochains and Homotopy Type}, Publ. Math. IHES, 103 (2006), 213-246.

\bibitem{Mane}M. Manetti, \textit{Extended deformation functors}, Int. Math. Res. Not. 14 (2002), 719-756.

\bibitem{Mar}M. Markl, \textit{Intrinsic brackets and the L-infinity deformation theory of bialgebras}, J. Homotopy
Relat. Struct. 5 (2010), no. 1, 177–212.

\bibitem{May} J. P. May, {\it The geometry of iterated loop spaces}, Springer, Berlin, 1972.

\bibitem{MV1}S. Merkulov, B. Vallette, \textit{Deformation theory of representation of prop(erad)s I}, J. für die reine und angewandte Math. (Crelles Journal), Issue 634 (2009), 51-106.

\bibitem{MV2}S. Merkulov, B. Vallette, \textit{Deformation theory of representation of prop(erad)s II}, J. für die reine und angewandte Math. (Crelles Journal), Issue 636 (2009), 125-174.

\bibitem{Mer1}S. Merkulov, \textit{Prop profile of Poisson geometry},  Comm. Math. Phys. 262 (2006), 117–135. 

\bibitem{Mer2}S. Merkulov, \textit{Formality theorem for quantization of Lie bialgebras},  Lett. Math.
 Phys. 106 (2016), no. 2, 169–195.
 
\bibitem{Mer3}S. Merkulov, \textit{Lectures on PROPs, Poisson geometry and deformation quantization}, Poisson geometry in mathematics and physics, 223–257,
Contemp. Math., 450, Amer. Math. Soc., Providence, RI, 2008.

\bibitem{Mer4}S. Merkulov, \textit{Graph complexes with loops and wheels}, Algebra, arithmetic, and geometry: in honor of Yu. I. Manin. Vol. II, 311–354,
Progr. Math., 270, Birkhäuser Boston, Inc., Boston, MA, 2009.

\bibitem{Mur}F. Muro, \textit{Moduli spaces of algebras over nonsymmetric operads}, Algebr. Geom. Topol. 14 (2014) 1489–1539.

\bibitem{PTVV}T. Pantev, B. Toen, M. Vaquié, G. Vezzosi, \textit{Shifted symplectic structures},  Inst. Hautes Études Sci. Publ. Math. 117 (2013), 271–328.

\bibitem{Preygel}{Pr} A. Preygel, \textit{Thom-Sebastiani and Duality for Matrix Factorizations, and Results on
the Higher Structures of the Hochschild Invariants}, Thesis (Ph.D.), M.I.T. 2012

\bibitem{Pri}J. P. Pridham, \textit{Unifying derived deformation theories}, Adv. Math. 224 (2010), 772–826.

\bibitem{Pri2}J. P. Pridham, \textit{Presenting higher stacks as simplicial schemes}, Adv. Math. 238 (2013), 184-245.

\bibitem{RT}N. Reshetikhin, V. G. Turaev, \textit{Invariants of 3-manifolds via link polynomials and quantum groups},  Invent. Math. 103 (1991), no. 3, 547–597.

\bibitem{Rez}C. W. Rezk, \textit{Spaces of algebra structures and cohomology of operads}, Thesis, MIT, 1996.

\bibitem{SW}P. Salvatore, N. Wahl, \textit{Framed discs operads and Batalin-Vilkovisky algebras}, Quart. J. Math. 54 (2003), 213-231.

\bibitem{Sev}P. Severa, \textit{Formality of the chain operad of framed little disks}, Lett. Math. Phys. 93 (2010), 29–35.

\bibitem{Sim}C. Simpson, \textit{Moduli of representations of the fundamental group of a smooth projective variety}, Publ. Math. IHES 80 (1994), 5-79.

\bibitem{Shoikhet} B. Shoikhet, \textit{Tetramodules over a bialgebra form a 2-fold monoidal category},
Appl. Categ. Structures {\bf 21} (2013), no.~3, 291--309.

\bibitem{Sul} D. Sullivan, \textit{Infinitesimal computations in topology}, Inst. Hautes Études Sci. Publ. Math. 47 (1977), 269–331.

\bibitem{Taillefer-HopfBimod} R. Taillefer, \textit{Cohomology theories of Hopf bimodules and cup-product},
Algebr. Represent. Theory {\bf 7} (2004), no.~5, 471--490.

\bibitem{Tam1}D. Tamarkin, \textit{Another proof of M. Kontsevich formality theorem}, preprint arXiv:math/9803025, 1998.

\bibitem{Tam2}D. Tamarkin. \textit{Formality of chain operad of little discs}, Lett. Math. Phys. 66 (2003), 65–72.

\bibitem{Ta-deformationofd-algebra}  D. Tamarkin, \textit{Deformation complex of a d-algebra is a (d+1)-algebra}, preprint arXiv:math/0010072.

\bibitem{TV}B. Toën, G. Vezzosi, \textit{Homotopical algebraic geometry II. Geometric stacks and applications},
Mem. Amer. Math. Soc. 193 (2008), no. 902, x+224 pp.

\bibitem{Toe}B. Toën, \textit{Derived algebraic geometry}, EMS Surv. Math. Sci. 1 (2014), no. 2, 153–240.

\bibitem{Toen-ICM} B. Toën, \textit{Derived Algebraic Geometry and Deformation Quantization}, ICM lecture.

\bibitem{To-Bourbaki}  B. . Toën, \textit{Probl\`emes de modules formels}, S\'eminaire BOURBAKI
68\`eme ann\'ee, 2015-2016, no 1111.

\bibitem{Tur} V. G. Turaev, \textit{Skein quantization of Poisson algebras of loops on surfaces}, Ann. Sci. École Norm. Sup. (4) 24 (1991), no. 6, 635–704.

\bibitem{Val}B. Vallette, \textit{A Koszul duality for props}, Trans. Amer. Math. Soc. 359 (2007), 4865-4943.

\bibitem{Wit}E. Witten, \textit{Quantum field theory and the Jones polynomial}, Commun. Math. Phys. 121 (1989), 351-399.

\bibitem{Yal2}S. Yalin, \textit{Classifying spaces of algebras over a prop},  Algebr. Geom. Topol. 14, issue 5 (2014), 2561--2593.

\bibitem{Yal3}S. Yalin, \textit{Simplicial localization of homotopy algebras over a prop},  Math. Proc. Cambridge Philos. Soc.   157 (2014), no. 3, 457–468.

\bibitem{Yal4}S. Yalin, \textit{Maurer-Cartan spaces of filtered $L_{\infty}$ algebras}, J. Homotopy Relat. Struct. 11 (2016), 375–407. 

\bibitem{Yal5}S. Yalin, \textit{Moduli stacks of algebraic structures and deformation theory}, J. Noncommut. Geom. 10 (2016), 579-661.

\bibitem{Yal6}S. Yalin, \textit{Function spaces and classifying spaces of algebras over a prop}, Algebr. Geom. Topology 16 (2016), 2715–2749.

\bibitem{Yal7}S. Yalin, \textit{Realization spaces of algebraic structures on chains}, to appear in Int. Math. Res. Not.

\bibitem{GY}G. Ginot, S. Yalin, \textit{Deformation theory of bialgebras, higher Hochschild cohomology and formality}, preprint arXiv:1606.01504.

\end{thebibliography}
\end{document}